\documentclass[12pt]{amsart}

\usepackage{amsmath}
\usepackage{amssymb}
\usepackage{amscd}
\usepackage{latexsym}
\usepackage{amsthm}
\usepackage{mathrsfs}
\usepackage{color}
\usepackage[backend=bibtex,style=alphabetic,isbn=false,url=false,maxnames=5,firstinits=true]{biblatex}
\usepackage{amsfonts}
\usepackage{enumitem}
\usepackage{array}
\usepackage{setspace}
\usepackage{tikz}
\usepackage{tikz-cd}
\usetikzlibrary{arrows}
\renewbibmacro{in:}{}

\usepackage{stmaryrd,graphicx}

\newtheorem{thm}{Theorem}[section]
\newtheorem*{thm*}{Theorem}

\newtheorem*{cor*}{Corollary}
\newtheorem{lem}[thm]{Lemma}
\newtheorem*{lem*}{Lemma}
\newtheorem{prop}[thm]{Proposition}
\newtheorem*{prop*}{Proposition}

\theoremstyle{definition}
\newtheorem{defn}{Definition}[section]
\newtheorem*{defn*}{Definition}

\theoremstyle{remark}
\newtheorem{rem}{Remark}[section]
\newtheorem*{rem*}{Remark}

\newtheorem*{problem*}{Problem}

\newcommand{\Q}{\mathbb Q}
\newcommand{\R}{\mathbb R}

\newcommand{\CC}{\mathscr C}
\newcommand{\B}{\mathscr B}

\newcommand{\WP}{\mathcal{WP}}
\newcommand{\Z}{\mathbb Z}

\newcommand{\bA}{\mathbb A}

\newcommand{\DD}{\mathbb D}
\newcommand{\TT}{\mathbb T}

\newcommand{\wt}{\widetilde}

\newcommand{\eps}{\epsilon}

\newcommand{\old}{\mathrm{\,old}}

\newcommand{\upto}{\nearrow}
\newcommand{\downto}{\searrow\,}

\DeclareMathOperator{\opnext}{next}
\DeclareMathOperator{\pExp}{Exp}

\DeclareMathOperator{\Sym}{Sym}
\DeclareMathOperator{\inv}{inv}
\DeclareMathOperator{\area}{area}
\DeclareMathOperator{\dinv}{dinv}
\DeclareMathOperator{\tdinv}{tdinv}
\DeclareMathOperator{\maxtdinv}{maxtdinv}

\DeclareMathOperator{\Id}{Id}

\title{Toric braids and $(m,n)$-parking functions}
\author{Anton Mellit}
\begin{document}

\onehalfspacing

\begin{abstract}
The Dyck path algebra construction from \cite{carlsson2015proof} is interpreted as
a representation of ``the positive part'' of the group of toric braids.
Then certain sums over $(m,n)$-parking functions are related to evaluations
 of this 
representation on some special braids. The compositional $(km,kn)$-shuffle
 conjecture of \cite{bergeron2015compositional} is then shown to be a corollary
  of this relation.
\end{abstract}

\maketitle

\section{Introduction}
in \cite{carlsson2015proof} Erik Carlsson and the author proved a generalization of the shuffle conjecture of \cite{haglund2005diagcoinv}. This conjecture gave a combinatorial interpretation of the Frobenius character of the algebra of diagonal coinvariants in $2$ sets of variables, in terms of parking functions, which are labellings of $(n,n)$-Dyck paths, i.e. lattice paths from $(0,0)$ to $(n,n)$ going only up and right and staying above the diagonal. In the papers \cite{gorsky2013compactified}, \cite{hikita2014affine}, \cite{gorsky2015refined} and \cite{bergeron2015compositional} this conjecture was further extended to include parking functions on $(m,n)$-Dyck paths, which are paths from $(0,0)$ to $(m,n)$ staying above the $(m,n)$-diagonal.

In an attempt to prove the $(m,n)$-shuffle conjecture the author discovered a precise relation between toric braids acting in a certain representation and parking functions (see Theorem \ref{thm:main}). For an $(m,n)$-Dyck path and a real number $h>0$ we slice the path by a line of slope $n/m$ at height $h$ and obtain a union of intervals $c$, which we call a coloring of the line. Projecting these intervals to the punctured torus $\TT_0 = (\R/\Z)^2\setminus\{(0,0)\}$ we obtain a braid $B$, i.e. an evolution of a tuple of distinct points on $\TT_0$. The construction of Dyck path algebra of \cite{carlsson2015proof} is interpreted as a representation of a certain submonoid of the group of braids on $\TT_0$. Then it is shown that the braid $B$, when evaluated in this representation, gives a certain sum over parking functions based on Dyck paths with slice $c$.

The full statement of the compositional $(m,n)$-shuffle conjecture of \cite{bergeron2015compositional} is then deduced from Theorem \ref{thm:main}.

Performing computations with braids we exploit two different approaches. One approach is to follow the braid ``top to bottom'' as the time varies. In this way we obtain the usual expressions of braids in terms of the generators of the braid group. Another, a completely orthogonal and not so standard approach is to move the braid sideways, following it as it breaks, or strands are created or deleted. As shown in this work, the second approach is the one that produces Dyck paths and parking functions. Perhaps this idea can be used in other situations.

The ``double Dyck path algebra'' $\bA_{q,t}$ that we use throughout the paper contains copies of the double affine Hecke algebra (DAHA) for $\operatorname{GL}_n$ for all $n$, related between themselves by raising and lowering operators. Existence of these operators is very useful for the present paper. No attempt was made to compare techniques of the present work with previous works on DAHA, i.e. \cite{gorsky2014torus} and \cite{cherednik2014DAHA}. We can expect that such a comparison would produce many interesting results.

I am grateful to Adriano Garsia for a suggestion to read about relationships between shuffle algebras and toric knots in \cite{morton2014HOMFLYPT}. This suggestion eventually lead me to a right guess for a generalization of the $(m,n)$-shuffle conjecture, which was general enough to satisfy recursions. I thank Shehryar Sikander for discussions about toric braid groups and  pointing out the reference \cite{bellingeri2004presentations}, and Fernando Rodriguez-Villegas for interesting discussions about the subject and continuous encouragement. Misha Mazin and Andrei Negut explained some of their work to me, which was helpful and interesting. I also thank Erik Carlsson for introducing me to shuffle conjectures, and for interesting discussions. The present research was performed during the author's stay at SISSA, and partly at ICTP in Trieste. I am thankful to these institutions for the support and stimulating environment.

\section{Notations}
Notations are similar to the ones in \cite{carlsson2015proof}. The base ring is $\Q(q,t)$. The ring of symmetric functions is denoted by $\operatorname{Sym}$ or $\operatorname{Sym}[X]$. For $F\in\Sym[X]$, $F$ and $F[X]$ is the same thing. Square brackets are used to denote plethystic substitution as in $F[(q-1)X]$. Plethystic exponential is defined as
\[
\pExp[X] = \sum_{n=0}^\infty h_n[X] = \exp\left(\sum_{n=1}^\infty \frac{p_n[X]}n \right).
\]
For $k=0,1,2,\ldots$ we set 
\[
V_k:=\Sym[X]\otimes \Q(q,t)[y_1, y_2,\ldots,y_k],\quad V_*:=\bigoplus_{k=0}^\infty V_k,
\]
$V_k$ is a $\lambda$-ring in the usual way, so plethystic operations are defined on it.

We will often have elements $T_1, T_2, \cdots, T_{k-1}$ satisfying the classical braid group relations. The following shorthands are useful for $i\leq j$:
\[
T_{i\upto j} = T_i T_{i+1} \cdots T_{j-1},\quad T_{j\downto i} = T_{j-1} T_{j-2} \cdots T_{i},
\]
\[
T_{i\upto j}^* = T_i^{-1} T_{i+1}^{-1} \cdots T_{j-1}^{-1},\quad T_{j\downto i}^* = T_{j-1}^{-1} T_{j-2}^{-1} \cdots T_{i}^{-1},
\]
\[
T_{i\upto i} = T_{i\downto i} = T_{i\upto i}^* = T_{i\downto i}^* = \Id,
\]
and we sometimes extend the notation to include the cases $i> j$ as follows: 
\[
T_{i\upto j}:=T_{i\downto j}^*,\quad T_{j\downto i}:=T_{j\upto i}^*\quad(i>j).
\]
For an expression $f(z)$ in the variable $z$ the following denotes the coefficient of $z^i$:
\[
f(z)\Big|_{z^i}.
\]

\section{Dyck path algebra}

\subsection{The algebra}
The Dyck path algebra was introduced in \cite{carlsson2015proof} to compute certain characteristic functions associated to Dyck paths. Here is the definition:
\begin{defn}[\cite{carlsson2015proof}]
The Dyck path algebra $\bA_{q}$ (over $\Q(q)$) is the path algebra of the quiver with vertex set $\Z_{\geq 0}$, arrows $d_+$ from $k$ to $k+1$, arrows $d_-$ from $k+1$ to $k$, and loops $T_1, T_2, \dots, T_{k-1}$ from $k$ to $k$ ($k\geq 0$) subject to the following relations:
\begin{equation}\label{eq:Trel}
(T_i-1)(T_i+q)=0, \quad T_i T_{i+1} T_i = T_{i+1} T_i T_{i+1}, \quad T_i T_j = T_j T_i \quad (|i-j|>1),
\end{equation}
\begin{equation}\label{eq:Tdminusrel}
d_-^2 T_{k-1} = d_-^2, \quad T_i d_- = d_- T_i \quad(1\leq i \leq k-1),
\end{equation}
\begin{equation}\label{eq:Tdplus}
T_1 d_+^2 = d_+^2,\quad d_+ T_i = T_{i+1} d_+ \quad(1\leq i \leq k),
\end{equation}
\begin{equation}\label{eq:extrarel1}
d_- (d_+ d_- - d_- d_+) T_{k-1} = q (d_+ d_- - d_- d_+) d_-\quad (k\geq 2),
\end{equation}
\begin{equation}\label{eq:extrarel2}
T_1 (d_+ d_- - d_- d_+) d_+ = q d_+ (d_+ d_- - d_- d_+) \quad (k\geq 1),
\end{equation}
where in each identity $k$ denotes the index of the vertex where the respective paths begin. The idempotent corresponding to the vertex $k$ will be denoted $e_k$.
\end{defn}

It turns out that $\bA_q$ contains a copy of the affine Hecke algebra of type $\widetilde{A_{n-1}}$ for each $n>0$ :

\begin{prop}[\cite{carlsson2015proof}]\label{lem:yfromd}
For each $k>0$ define loops $y_1, y_2, \ldots, y_k$ from $k$ to $k$ by the relations 
\begin{equation}\label{eq:rely1}
y_1 = \frac{1}{q^{k-1}(q-1)} (d_+ d_- - d_- d_+)\; T_{k\downto 1} ,
\end{equation}
\begin{equation}\label{eq:rely2}
y_{i+1} = q T_i^{-1} y_{i} T_i^{-1} \quad (1\leq i\leq k-1).
\end{equation}
Then the following identities hold:
\begin{align}\label{eq:relTy}
y_i T_j =& T_j y_i  &(i\notin \{j, j+1\}),
\\ \label{eq:relydminus}
y_i d_- =& d_- y_i  &(1\leq i \leq k-1),
\\ \label{eq:relydplus}
d_+ y_i =& T_{1\upto i+1}\; y_i\; T_{i+1\downto 1}^* d_+,
\qquad&(1\leq i \leq k),
\\ \label{eq:relycomm}
y_i y_j =& y_j y_i &(1\leq i,j \leq k).
\end{align}
\end{prop}

\subsection{Two actions}
In \cite{carlsson2015proof} actions of $\bA_q$ and $\bA_{q^{-1}}$ were constructed on the sequence of spaces $V_* = \bigoplus_{k=0}^\infty V_k$, where $V_k=\Sym[X]\otimes \Q(q, t) [y_1, \ldots, y_k]$\footnote{The first action does not use the variable $t$, so in fact is an action on $\Sym[X]\otimes \Q(q) [y_1, \ldots, y_k]$.} For the purpose of this paper it is convenient to consider slightly different actions. We first formulate our statements in terms of the modified actions, and then  explain how to obtain the modified results from the ones of \cite{carlsson2015proof}.

\begin{prop}\label{lem:action1}
The following operations define an action of $\bA_q$ on $V_*$:
\[
T_i F = \frac{(q-1) y_i F + (y_{i+1} - q y_i) \operatorname{swap}_{y_i\leftrightarrow y_{i+1}} F}{y_{i+1} - y_i},
\]
\[
d_- F = F[X- (q-1) y_k;y_1,\ldots,y_k] \pExp[-y_k^{-1} X]\Big|_{y_k^0},
\]
\[
d_+ F = -T_{1\upto k+1}\; (y_{k+1} F[X+(q-1) y_{k+1}; y_1, \ldots, y_k]).
\]
Here $\operatorname{swap}_{y_i\leftrightarrow y_{i+1}}$ is the operator of interchanging the variables $y_i$ and $y_{i+1}$. Moreover, the operators $y_i\in \bA_q$ act precisely by multiplication by $y_i$.
\end{prop}

\begin{prop}\label{lem:action2}
Define $d_+^*:V_k \to V_{k+1}$ by 
\[
d_+^* F = \gamma\left(F[X+(q-1) y_{k+1}; y_1, \ldots, y_k])\right),
\]
where $\gamma$ is the operator which sends $y_i$ to $y_{i+1}$ for $i<k+1$ and $y_{k+1}$ to $t y_1$. Then the operators $T_i^{-1}, d_-, d_+^*$ define an action of the conjugate algebra $\bA_{q^{-1}}$ on $V_*$.
\end{prop}

The action of the elements $y_i\in \bA_{q^{-1}}$ defines commuting operators $z_i$. Specifically, we have
\[
z_1 = \frac{q^k}{(1-q)} (d_+^* d_- - d_- d_+^*)\; T_{k\downto 1}^*,
\]
\[
z_{i+1} = q^{-1} T_i z_i T_i \quad (1\leq i\leq k-1).
\]

Speaking about actions of $\bA_q$ on $V_*$ we always consider only actions such that $e_k$ is the projection on $V_k$.
Let $\rho, \rho^*$ be actions of $\bA_q$, $\bA_{q^{-1}}$ on $V_*$ which respect the decomposition $V_*$.
We denote $\rho(T_i)$, $\rho(d_-)$, $\rho(d_+)$, $\rho(y_i)$, $\rho^*(d_+)$, $\rho^*(y_i)$ simply by $T_i$, $d_-$, $d_+$, $y_i$, $d_+^*$, $z_i$ respectively.
\begin{defn}\label{defn:correctinter}
We say that $\rho$, $\rho^*$ are \emph{correctly intertwined} if the following identities hold for all $k\geq 0$. 
\[
\rho^*(T_i) = T_i^{-1} \quad(1\leq i\leq k-1),\quad \rho^*(d_-) = d_-,
\]
\[
d_+ z_i = z_{i+1} d_+, \quad d_+^* y_i = y_{i+1} d_+^* \quad(1\leq i \leq k),
\]
\[
z_1 d_+ = - t q^{k+1} y_1 d_+^* : V_k\to V_{k+1}.
\]
\end{defn}

Then we have
\begin{prop}\label{lem:intertwined}
The actions of Propositions \ref{lem:action1}, \ref{lem:action2} are correctly intertwined.
\end{prop}

This can be reformulated as follows. Let $\bA_{q,t}$ be the path algebra of the quiver with vertex set $\Z_{\geq 0}$, arrows $d_+$ and $d_+^*$ from $k$ to $k+1$, arrows $d_-$ from $k+1$ to $k$, and loops $T_1, T_2, \dots, T_{k-1}$ from $k$ to $k$ ($k\geq 0$) subject to the following relations:
\begin{enumerate}
\item relations of $\bA_q$ for $T_i, d_-, d_+$,
\item relations of $\bA_{q^{-1}}$ for $T_i^{-1}, d_-, d_+^*$,
\item relations of Definition \ref{defn:correctinter}.
\end{enumerate}
Then the actions of Propositions \ref{lem:action1} and \ref{lem:action2} form an action of $\bA_{q,t}$.

\subsection{Proofs of the modified relations}
The following relations do not require any adaptation: \eqref{eq:Trel}, the second parts of \eqref{eq:Tdminusrel} and \eqref{eq:Tdplus}, \eqref{eq:rely2}\footnote{Here $y_i$ are the operators of multiplication by $y_i$, so before we prove Proposition \ref{lem:action1} we cannot simply deduce these relations from Proposition \ref{lem:yfromd}} \eqref{eq:relTy}, \eqref{eq:relydminus}, \eqref{eq:relydplus} and \eqref{eq:relycomm}.
For the purpose of the proof denote the operators $d_-$, $d_+$,\dots from \cite{carlsson2015proof} by $d_-^\old$, $d_+^\old$ and so on. We clearly have 
\[
T_i^\old = T_i, \quad y_i^\old = y_i, \quad d_+^{*\old}=d_+^*, \quad d_-^\old = - d_- y_k,
\]
and
\[
d_+ = -T_{1\upto k+1}\; y_{k+1}\; T_{k+1\downto 1}^*\, d_+^\old.
\]
The strategy is to express everything in terms of $d_-$ and $d_+^\old$.
Let $\tau$ be the shift operator acting on the $\Sym[X]$ part of $V_k$: $F[X]\to F[X+1]$. Then $\tau$ commutes with $d_+$, $d_+^\old$, $d_+^*$ , $T_i$ and $y_i$. We have
\[
[\tau, d_-^\old] = d_- \tau.
\]
Therefore applying $[\tau, -]$ to the old relation
\[
y_1 = -\frac{1}{q^{k-1}(q-1)} (d_-^\old d_+^\old - d_+^\old d_-^\old)\, T_{k \downto 1}
\]
we obtain
\[
0 = -\frac{1}{q^{k-1}(q-1)} (d_- d_+^\old - d_+^\old d_-)\, T_{k\downto 1}\, \tau,
\]
so $d_-$ commutes with $d_+^\old$:
\[
d_- d_+^\old = d_+^\old d_-.
\]
Moreover we easily have from the definition
\[
d_-\, T_{k+1\downto 1}^*\, d_+^\old = \Id.
\]
Therefore we can express
\[
d_+ d_- - d_- d_+ = - \;T_{1\upto k}\, y_k \,T_{k\downto 1}^*\, d_+^\old d_- \,+\, d_- \,T_{1\upto k+1} \,y_{k+1}\, T_{k+1\downto 1}^*\, d_+^\old
\]
\[
=T_{1\upto k}\, y_k d_- \left(- T_{k\downto 1}^* \,+\, q T_{k}^{-2} \, T_{k\downto 1}^* \right) d_+^\old
\;=\; (q-1) T_{1\upto k} y_k.
\]
This establishes \eqref{eq:rely1}. Now it is straightforward to use it to prove \eqref{eq:extrarel1} and \eqref{eq:extrarel2}. 

The first parts of \eqref{eq:Tdminusrel} and \eqref{eq:Tdplus} can be proved similarly to the proof in \cite{carlsson2015proof}. At this point Proposition \ref{lem:action1} is proved.
\begin{rem}\label{rem:oldnew}
Let $M:V_k\to V_k$ be the operator of multiplication by $(-1)^k y_1 y_2\cdots y_k$ for each $k\geq 0$. Then we have
\[
M d_-^\old = d_- M, \quad M d_+^\old = d_+ M, \quad M y_i^\old = y_i M
\]
easily from the definitions. If $M$ was surjective, we could deduce our statements from the old ones, but $M$ is not surjective. However, it is injective, so the old statements can be deduced from the new ones.
\end{rem}

For Proposition \ref{lem:action2} we only need to show that the relations \eqref{eq:extrarel1} and \eqref{eq:extrarel2} with $q$ replaced by $q^{-1}$ hold for $d_-$ and $d_+^*$. The latter can be obtained by an application of $[\tau, -]$ to the corresponding old relation. To prove \eqref{eq:extrarel1} we follow the proof in \cite{carlsson2015proof} and replace each occurrence of the operator $B_i$ by $-B_{i-1}$.

Finally, we prove Proposition \ref{lem:intertwined}. Essentially there are three statements:
\[
d_+^* y_i = y_{i+1} d_+^*
\]
is the same as in \cite{carlsson2015proof}. The other two are about $z_i$ and $d_+$. First we relate $z_1$ and $z_1^\old$:
\[
z_1^\old = \frac{q^k}{1-q} (-d_+^* d_- y_k + d_- y_{k+1} d_+^*)\,T_{k\downto  1}^*\,=\, -z_1 \,T_{1\upto k}\, y_k T_{k\downto 1}^*.
\]
Thus we have
\[
z_1 d_+ = z_1^\old d_+^\old = -t q^{k+1} y_1 d_+^*.
\]
To prove $d_+ z_i = z_{i+1} d_+$ recall from \cite{carlsson2015proof} that this is equivalent to 
\[
T_1^{-1} d_+(d_+^* d_- - d_- d_+^*) = (d_+^* d_- - d_- d_+^*) d_+.
\]
Applying $[\tau, -]$ to the corresponding old relation gives
\[
T_1^{-1} d_+^\old(d_+^* d_- - d_- d_+^*) = (d_+^* d_- - d_- d_+^*) d_+^\old.
\]
Now we multiply both sides on the left by $T_2 \cdots T_k y_{k+1} T_k^{-1} \cdots T_2^{-1}$ and obtain the desired relation. This completes the proof of Proposition \ref{lem:intertwined}.

\subsection{Structure theorems}
In \cite{carlsson2015proof} the following theorem was established for the algebra $\bA_q$:
\begin{thm}[\cite{carlsson2015proof}]
As an $\bA_q$-module, the spherical module $\bA_q e_0$ is isomorphic to $V_*$ with the old $\bA_q$-action. For the new action the statement is
\[
\bA_q e_0 \cong \bigoplus_{k=0}^\infty y_1 \cdots y_k V_k.
\]
The spherical algebra $e_0 \bA_q e_0$ is isomorphic to the commutative algebra of symmetric functions over $\Q(q)$, 
\[
e_0 \bA_q e_0 \cong \Sym[X].
\]
\end{thm}
The new statement follows from the old one by Remark \ref{rem:oldnew}.

For the algebra $\bA_{q,t}$ the corresponding statement is a little more complicated than the one in \cite{carlsson2015proof}
\begin{thm}\label{thm:structure2}
The kernel of the natural homomorphism of $\bA_{q,t}$-modules $\bA_{q,t} e_0 \to V_*$ is
\[
I := \langle (d_- d_+^* - 1) d_+^{*k} e_0, \; (d_+ + q^k y_1 d_+^*) d_+^{*k} e_0 \;|\;k\geq 0\rangle,
\]
and the resulting map $\bA_{q,t} e_0/I \to V_*$ is an isomorphism.
\end{thm}
\begin{proof}
Denote $\bA_{q,t} e_0/I \to V_*$ by $\varphi$.
The strategy, as in  \cite{carlsson2015proof} is to construct a map $\psi: V_* \to \bA_{q,t} e_0$ such that $\varphi \psi=\Id$, and then prove that $\psi$ is surjective by induction. By Lemma 5.5 in \cite{carlsson2015proof} elements of the following form form a basis of $V_k$:
\[
\varphi(d_-^m y_1^{a_1} \cdots y_{k+m}^{a_{k+m}} d_+^{* k+m} e_0)\quad (m\geq 0, a_i\geq 0, a_{k+1}\geq \cdots\geq a_{k+m}>0).
\]
Note that we had to translate the statement to the new generators and use the fact that $d_+^{*k+m}=1\in V_{k+m}$. Thus we can define $\psi$ in such a way that 
\begin{equation}\label{eq:basis}
\psi \varphi(d_-^m y_1^{a_1} \cdots y_{k+m}^{a_{k+m}} d_+^{* k+m} e_0) = 
d_-^m y_1^{a_1} \cdots y_{k+m}^{a_{k+m}} d_+^{* k+m} e_0.
\end{equation}
It remains to show that the image of $\psi$ spans $\bA_{q,t}$ modulo $I$.

Similarly to the proof of Theorem 7.1 of \cite{carlsson2015proof} we show that applying $d_+^*$ and $d_+$ to the elements of the form \eqref{eq:basis} we can express the result in the same form. The only difference is that applying $d_-$ can make $a_{k+m}=0$. In this case we can simply decrease $m$ by $1$ using the first set of generators of $I$.
\end{proof}

\subsection{Characteristic functions}\label{sec:charfunc}
The isomorphism $e_0 \bA_q e_0 \cong \Sym[X]$ possesses a combinatorial interpretation. Let $n\geq 0$, and let $\pi$ be an $(n,n)$-Dyck path, i.e. a lattice path from $(0,0)$ to $(n, n)$ consisting of North and East steps, which stays weakly above the diagonal. For the purposes of this paper it is useful to endow Dyck paths with an extra structure. A \emph{marking} of a Dyck path $\pi$ is a choice of a subset $S$ of the set of \emph{corners} of $\pi$, i.e. cells $(i,j)$ which are above the path, but $(i, j-1)$ and $(i+1, j)$ are both below the path. For a Dyck path $\pi$ and a marking $S$ we define the \emph{characteristic function} as follows:
\[
\chi(\pi, S):=\sum_{w\in \Z_{>0}^n \;\text{$S$-admissible}} q^{\inv(\pi, w)} x_{w_1} x_{w_2} \cdots x_{w_n} \in \Sym[X]\otimes\Q(q),
\]
where $S$-admissible means $w_i>w_j$ for each $(i,j)\in S$ and 
\[
\inv(\pi, w) = \#\{(i,j): w_i>w_j \;\text{and $(i,j)$ is weakly below $\pi$}\}.
\]
This defines a symmetric function of degree $n$, and the desired combinatorial interpretation is:

\begin{thm}[\cite{carlsson2015proof}]\label{thm:combinatorial}
Via the isomorphism $e_0 \bA_q e_0 \cong \Sym[X]$, for any marked Dyck path $(\pi, S)$, $\chi(\pi, S)$ corresponds to an element of $e_0 \bA_q e_0$ constructed as follows. We follow the path from right to left applying
\begin{enumerate}
\item $\frac{1}{q-1}(d_- d_+ - d_+ d_-)$ for each marked corner,
\item $d_-$ for each vertical step which is not a part of a marked corner,
\item $d_+$ for each horizontal step which is not a part of a marked corner.
\end{enumerate}
\end{thm}

\subsection{Replication of actions}
It turns out that starting with the two actions constructed in Propositions \ref{lem:action1} and \ref{lem:action2} we can construct infinitely many actions using the following construction:
\begin{prop}\label{prop:replication}
Let $\rho$, $\rho^*$ be actions of $\bA_q$, $\bA_{q^{-1}}$ resp. which are correctly intertwined. Denote $\rho(T_i)$, $\rho(d_-)$, $\rho(d_+)$, $\rho(y_i)$, $\rho^*(d_+)$, $\rho^*(y_i)$ by $T_i$, $d_-$, $d_+$, $y_i$, $d_+^*$, $z_i$ resp. Then the assignments $\rho_1$, $\rho_2$
\[
\rho_1(d_-) = \rho_2(d_-)=d_-,\quad \rho_1(T_i)=T_i, \quad \rho_2(T_i)=T_i^{-1},
\]
\[
\rho_1(d_+) = -(qt)^{-1} z_1 d_+,\quad \rho_2(d_+) = -y_1 d_+^*
\]
induce actions of $\bA_q$, $\bA_{q^{-1}}$ resp. Moreover, the pairs $(\rho, \rho_2)$ and $(\rho_1, \rho^*)$ are correctly intertwined. \footnote{The Proposition remains true if we multiply $\rho_1(d_+)$, $\rho_2(d_+)$ by arbitrary constants.}
\end{prop}
\begin{proof}
From the obvious symmetry of the statement it is enough to check only half. So we verify that $\rho_1$ induces an action of $\bA_q$. The equations \eqref{eq:Trel}, \eqref{eq:Tdminusrel} do not need verification. The equations \eqref{eq:Tdplus} reduce to
\[
T_1 z_1 d_+ z_1 d_+ = T_1 z_1 z_2 d_+^2 = z_1 z_2 d_+^2 = z_1 d_+ z_1 d_+,
\]
\[
z_1 d_+ T_i = z_1 T_{i+1} d_+ = T_{i+1} z_1 d_+.
\]
Equations \eqref{eq:extrarel1}, \eqref{eq:extrarel2} are also verified in the same way. Thus the actions are defined correctly. Notice that the same approach shows
\[
\rho_1(y_1) = z_1 y_1.
\]

Next we check the conditions of Definition \ref{defn:correctinter}. These amount to
\[
z_1 d_+ z_i = z_{i+1} z_1 d_+, \quad d_+^* z_1 y_1 = q T_1^{-1} z_1 y_1 T_1^{-1} d_+^*,
\]
\[
z_1 z_1 d_+ = -t q^{k+1} z_1 y_1 d_+^*,
\]
which are also clear.
\end{proof}

Iterating the construction above we obtain
\begin{thm}\label{thm:actions}
For every $m,n\in \Z_{\geq0}$ such that $\operatorname{gcd}(m,n)=1$ we have actions $\rho_{m,n}$, $\rho_{m,n}^*$ of $\bA_q$, $\bA_{q^{-1}}$ respectively on $V_*$ so that:
\begin{itemize}
\item $\rho_{0,1}$ is the action of Proposition \ref{lem:action1},
\item $\rho_{1,0}^*$ is the action of Proposition \ref{lem:action2},
\item For each $m,n$ we have $\rho_{m,n}(d_+) = -q^k \rho_{m,n}^*(d_+^*)$,
\item For each $m,n$ and $m', n'$ satisfying condition $m' n - m n' = 1$ the actions $\rho_{m,n}$ and $\rho_{m', n'}^*$ are correctly intertwined, and the actions $\rho_{m+m', n+n'}$, $\rho_{m+m', n+n'}^*$ are obtained from these by the construction of Proposition \ref{prop:replication}.
\end{itemize}
\end{thm}

\begin{figure}
\begin{tikzpicture}[scale=1.5]
\begin{scope}<+->;
  \draw[gray,step=1cm,very thin] (0,0) grid (5.5,5.5);

  \draw[gray,thick,->] (0, 0) -- (5.5, 0) node[black,right] {$m$};
  \draw[gray,thick,->] (0, 0) -- (0, 5.5) node[black,above] {$n$};

  \draw (1,1) node[draw,circle,inner sep=1pt,fill] {} node[above right] {$\rho_{1,1}$};
  \draw (1,0) node[draw,circle,inner sep=1pt,fill] {} node[above right] {$\rho_{1,0}$};
  \draw (0,1) node[draw,circle,inner sep=1pt,fill] {} node[above right] {$\rho_{0,1}$};
  \draw (2,1) node[draw,circle,inner sep=1pt,fill] {} node[above right] {$\rho_{2,1}$};
  \draw (1,2) node[draw,circle,inner sep=1pt,fill] {} node[above right] {$\rho_{1,2}$};
  \draw (1,3) node[draw,circle,inner sep=1pt,fill] {} node[above right] {$\rho_{1,3}$};
  \draw (2,3) node[draw,circle,inner sep=1pt,fill] {} node[above right] {$\rho_{2,3}$};
  \draw (1,4) node[draw,circle,inner sep=1pt,fill] {} node[above right] {$\rho_{1,4}$};
  \draw (3,1) node[draw,circle,inner sep=1pt,fill] {} node[above right] {$\rho_{3,1}$};
  \draw (3,2) node[draw,circle,inner sep=1pt,fill] {} node[above right] {$\rho_{3,2}$};
  \draw (2,5) node[draw,circle,inner sep=1pt,fill] {} node[above right] {$\rho_{2,5}$};
  \draw (3,5) node[draw,circle,inner sep=1pt,fill] {} node[above right] {$\rho_{3,5}$};
  \draw (3,4) node[draw,circle,inner sep=1pt,fill] {} node[above right] {$\rho_{3,4}$};
  \draw (4,1) node[draw,circle,inner sep=1pt,fill] {} node[above right] {$\rho_{4,1}$};
  \draw (5,2) node[draw,circle,inner sep=1pt,fill] {} node[above right] {$\rho_{5,2}$};
  \draw (5,3) node[draw,circle,inner sep=1pt,fill] {} node[above right] {$\rho_{5,3}$};
  \draw (4,3) node[draw,circle,inner sep=1pt,fill] {} node[above right] {$\rho_{4,3}$};
  \draw[thin,->] (0,0) -- (1,2);
  \draw[thin,->] (0,0) -- (2,3);
  \draw[thin,->] (0,0) -- (3,5);
  \draw[thin,->] (0,0) -- (1,3);
  \draw[thin,->] (0,0) -- (1,4);
  \draw[thin,->] (0,0) -- (2,5);
  \draw[thin,->] (0,0) -- (3,4);
  \draw[thin,->] (0,0) -- (1,1);
  \draw[thin,->] (0,0) -- (2,1);
  \draw[thin,->] (0,0) -- (3,2);
  \draw[thin,->] (0,0) -- (3,1);
  \draw[thin,->] (0,0) -- (4,1);
  \draw[thin,->] (0,0) -- (5,2);
  \draw[thin,->] (0,0) -- (5,3);
  \draw[thin,->] (0,0) -- (4,3);
\end{scope}
\end{tikzpicture}
\caption{The actions generated after $4$ iterations}
\label{fig:replication}
\end{figure}
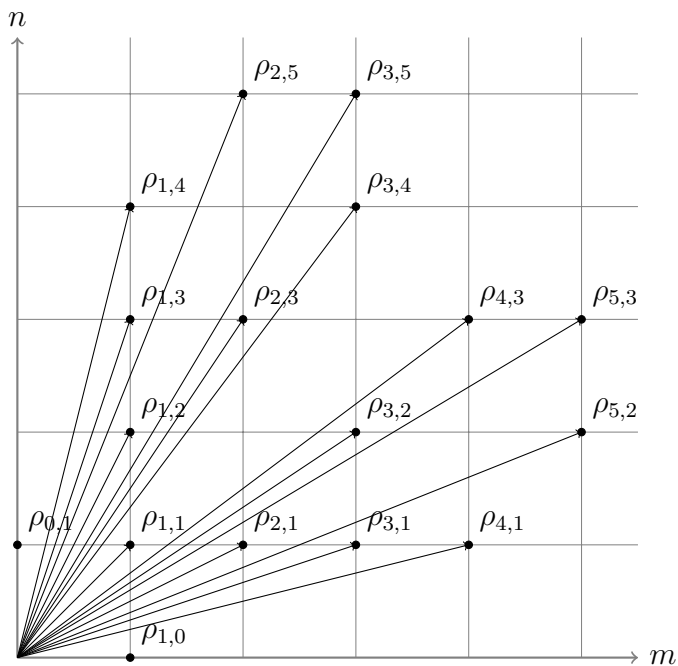

The process is demonstrated on Figure \ref{fig:replication}. Each time a new vector divides the sector between two neighboring vectors into two new sectors. After $4$ iterations we obtain the vectors shown on the picture. This is related to the Farey sequences. The fifth Farey sequence is
\[
0,\; \tfrac14,\; \tfrac13,\; \tfrac25,\; \tfrac12,\; \tfrac35,\; \tfrac23,\; \tfrac34,\; \tfrac11,\; \tfrac43,\; \tfrac32,\; \tfrac53,\; \tfrac21,\; \tfrac52,\; \tfrac31,\; \tfrac41,\; \tfrac10.
\]
The elements of the sequence are precisely the ratios $\frac{n}{m}$.

Another way to formulate Proposition \ref{prop:replication} is to say that we have two algebra homomorphisms $N, S:\bA_{q,t}\to \bA_{q,t}$ corresponding to the matrices (see \cite{bergeron2015compositional})
\[
N := \begin{pmatrix} 1& 1\\0 & 1\end{pmatrix}, \quad S:=\begin{pmatrix} 1 & 0\\ 1 & 1 \end{pmatrix},
\]
so that $N$ fixes the $d_-, d_+^*$-part and $S$ fixes the $d_-, d_+$-part. These matrices freely generate the monoid $\operatorname{SL}_2(\Z)_+\subset \operatorname{SL}_2(\Z)$ of matrices with non-negative entries (extended Euclidean algorithm), so we obtain an action of $\operatorname{SL}_2(\Z)_+$ on $\bA_{q,t}$, hence on its spherical part $e_0 \bA_{q,t} e_0$, which acts on $\Sym[X]$ by certain operators \footnote{In \cite{bergeron2015compositional} it is claimed that the whole group $\operatorname{SL}_2(\Z)$ acts on operators. This claim is false.}.

\subsection{Relation with $N$ and $S$ operators}
Next we want to show that the action of $N$ and $S$ on the spherical part $e_0 \bA_{q,t} e_0$ matches the one from \cite{bergeron2015compositional} and \cite{bergeron2014some}. First we treat the operator $N$.
\begin{prop}
Let $\nabla$ be the operator on $\Sym[X]$ defined in the Macdonald basis by
\[
\nabla \wt{H}_\lambda = (-1)^{|\lambda|} q^{n'(\lambda)} t^{n(\lambda)} \wt{H}_\lambda.
\]
Then for any $L\in e_0 \bA_{q,t} e_0$ we have the following identity between operators acting on $\Sym[X]$:
\[
\nabla L \nabla^{-1} = N(L).
\]
\end{prop}
\begin{proof}
Notice that the homomorphism $\bA_{q,t}e_0\to\bA_{q,t} e_0$ induced by $N$ preserves the module $I$ of Theorem \ref{thm:structure2}. Thus it induces an endomorphism of $V_k$, which we denote by $\nabla'$, such that
\[
\nabla' L = N(L) \nabla' \quad \text{for all $L\in\bA_{q,t}$}.
\]
It remains to show that $\nabla'=\nabla$ on $V_0$. It is enough to verify this identity on characteristic functions of $(n,n)$- Dyck paths for all $n$, because these span $\Sym[X]$. Let $\pi$ be such a Dyck path. By Theorem \ref{thm:combinatorial}
\[
\chi(\pi) = \pi(d_-^\old, d_+^\old) (1),
\]
where $\pi(d_-^\old, d_+^\old)$ means ``follow the path from right to left and apply $d_+^\old, d_-^\old$ for each horizontal resp. vertical step''. By Proposition 3.3 of \cite{carlsson2015proof} we have
\[
\bar\omega \chi(\pi) = \pi(d_-^\old, - q^{-k} d_+^\old) (1),
\]
where $-q^{-k} d_+^\old$ denotes the operator $V_*\to V_*$ which acts on $V_k$ as $-q^{-k} d_+^\old$. By Theorem 7.2 of \cite{carlsson2015proof} we have
\[
\nabla\chi(\pi) = \mathcal{N}\bar\omega\chi(\pi) =\pi(d_-^\old, -q^k d_+^*) (1).
\]
Conjugating by $M$ from Remark \ref{rem:oldnew} we obtain
\[
\nabla \chi(\pi) = \pi(d_-, q^k y_1 d_+^*).
\]
This equals to 
\[
\pi(d_-, -(qt)^{-1} z_1 d_+) = \nabla' \pi(d_-, d_+) = \nabla' \chi(\pi).
\]
\end{proof}

Let $\tau$ be the shift operator $F[X] \to F[X+1]$ on $\Sym[X]$ and $\tau^*$ be the operator of multiplication by $\pExp[-\frac{X}{(q-1)(t-1)}]$, which is conjugate to $\tau$ with respect to an appropriate modified Hall scalar product. Strictly speaking, $\tau^*$ does not act on $\Sym[X]$. It is an operator which maps $\Sym[X]$ to the space $\Sym[[X]]$, which is the completion of $\Sym[X]$ with respect to the grading by the degree. It is important to notice that $\tau^*$, and therefore the composition $\tau^* \tau$, is injective. Hence the following is true:

\begin{prop}
Let $D$ be a homogeneous operator acting on $\Sym[X]$. Then there exists at most one operator $S(D)$ on $\Sym[X]$ such that
\[
\tau^* \tau S(D) = D \tau^* \tau.
\]
Here we used a natural extension of $D$ to the space $\Sym[[X]]$.
Furthermore, operators $D$ such that $S(D)$ exists form a graded subalgebra of the algebra of all homogeneous operators, and $S$ is an algebra homomorphism.
\end{prop}

By a direct calculation we verify the following property of $S$:
\begin{prop}
As in \cite{bergeron2015compositional} define for any $n\in\Z$ an operator $D_n$ by the identity
\[
D_n F = F\left[X + (q-1)(t-1)z^{-1}\right] \pExp[-z X] \Big|_{z^n}.
\]
Then $S(D_n)$ exists and we have $S(D_n) = -D_{n+1}$.
\end{prop}

Thus we see that the operator $S$ from \cite{bergeron2015compositional} differs from ours only by a sign.

It remains to check that the endomorphism $S:\bA_{q,t}\to\bA_{q,t}$ indeed induces $S$ on an appropriate part of $e_0 \bA_{q,t} e_0$. Denote by $\bA_{q,t}^{<\infty}$ the subalgebra of $\bA_{q,t}$ generated by $T_i$, $d_-$, $d_+^*$ and $y_i$. It is clear that $\bA_{q,t}$ contains all the copies of $\bA_{q}$ constructed in Theorem \ref{thm:actions} with $(m,n)$ different from $(0,1)$.

\begin{prop}
For any element $L\in e_0 \bA_{q,t}^{<\infty} e_0$ we have the following identity on $\Sym[X]$:
\[
\tau^* \tau S(L) = L \tau^* \tau.
\]
\end{prop}
\begin{proof}
We extend the definition of $\tau^*$, $\tau$ to the spaces $V_k$ in the following way:
\[
\tau_k F = F[X+1] \prod_{i=1}^k (1-y_i^{-1}),
\]
\[
\tau_k^* F = F[X] \pExp\left[-\frac{X}{(q-1)(t-1)} - \frac{\sum_{i=1}^k y_i}{t-1}\right].
\]
In fact, $\tau_k$ is only defined on the subspace $y_1 y_2 \cdots y_k V_k \subset V_k$, which is enough for this proof.
These operators commute with $T_i$, $d_-$, $y_i$. As for $d_+^*$, we have
\[
d_+^* \tau_k^* = (1-y_1) \tau_{k+1}^* d_+^*,
\]
\[
(1- y_1^{-1}) d_+^* \tau_k = \tau_{k+1} d_+^*,
\]
therefore 
\[
\tau_{k+1}^* \tau_{k+1}  S(d_+^*)= -y_1\tau_{k+1}^* \tau_{k+1}  d_+^* = \tau_{k+1}^* (1 - y_1) d_+^* \tau_k = d_+^* \tau_k^* \tau_k.
\]
\end{proof}

The algebras $\bA_{q,t}$ and $e_0 \bA_{q,t} e_0$ have a bi-grading with the bi-degree of $d_+$, $d_+^*$ being $(0,1)$, $(1,0)$ respectively. The $\operatorname{SL}_2(\Z)_+$ permutes the graded components in the obvious way. Our action of $N$ sends the $(m,n)$ component to the $(m+n, n)$ component. The sign differs from the one of 
\cite{bergeron2015compositional} by $(-1)^n$. The action of $S$ sends $(m, n)$ to $(m, m+n)$ and the sign difference is $(-1)^m$. Thus multiplication by $(-1)^{m+n}$ intertwines our action and the one of \cite{bergeron2015compositional}.

Now we recall the operators 
\[
(C_a F)[X] = (-q)^{1-a} F[X + (q^{-1}-1)z] \pExp[z^{-1} X] z^a \Big|_{z^0}
\]
and the notation
\[
C_\alpha = C_{\alpha_1} C_{\alpha_2} \cdots C_{\alpha_r} 1 \quad(\alpha=(\alpha_1, \alpha_2,\ldots,\alpha_r)).
\]
We can write $C_\alpha$ with $\sum_{i=1}^r \alpha_i=k$ using our operators as follows:
\[
C_\alpha = (-q)^{r-k} \bar\omega( (-1)^r d_-^r y_1^{\alpha_1-1} \cdots y_r^{\alpha_r-1} d_+^r 1)
\]
\[
= (-1)^k q^{r-k} \rho_{0,1}^*(d_-^r y_1^{\alpha_1-1} \cdots y_r^{\alpha_r-1} d_+^r) 1.
\]
The left hand side of Conjecture 3.3 of \cite{bergeron2015compositional} becomes
\begin{equation}\label{eq:compshufflhs}
(-1)^{k(m+1)} q^{r-k} \rho_{m,n}^* (d_-^r y_1^{\alpha_1-1} \cdots y_r^{\alpha_r-1} d_+^r) 1.
\end{equation}
In the last equation the sign $(-1)^{k(n+1)}$ of \cite{bergeron2015compositional} was multiplied by $(-1)^{k(m+n)}$ and turned into $(-1)^{k(m+1)}$.
\begin{rem}
Working out the right hand side of \eqref{eq:compshufflhs} explicitly in terms of $d_+^*$ and $d_-$ produces constant term formulas as given in \cite{gorsky2015refined} and conjectured in \cite{bergeron2014some}, \cite{bergeron2015compositional}.
\end{rem}

\section{$(m,n)$-parking functions}
Fix relatively prime $m_1,n_1\in\Z_{>0}$. Denote the slope by $s:=\frac{n_1}{m_1}$. It is convenient to change the slope to $s_-:=s - \varepsilon$, where $\varepsilon$ is a very small positive constant we introduce solely for the purpose of breaking ties in comparisons between rational numbers and $s$. For instance, when working with diagrams contained in a very big square $N\times N$ it is enough to choose $\varepsilon = \frac{1}{2N^2}$, so that comparisons between different rational numbers whose numerator and denominator do not exceed $N$ are unaffected.

Let $g \in\Z_{>0}$ and put $m=g m_1$, $n=g n_1$. An $(m,n)$-Dyck path is a lattice path from $(0,0)$ to $(m,n)$ consisting of North and East steps  and staying weakly above the diagonal (equivalently, staying strictly above the line $y=s_- x$). The set of $(m,n)$-Dyck paths will be denoted $\DD_{m,n}$.

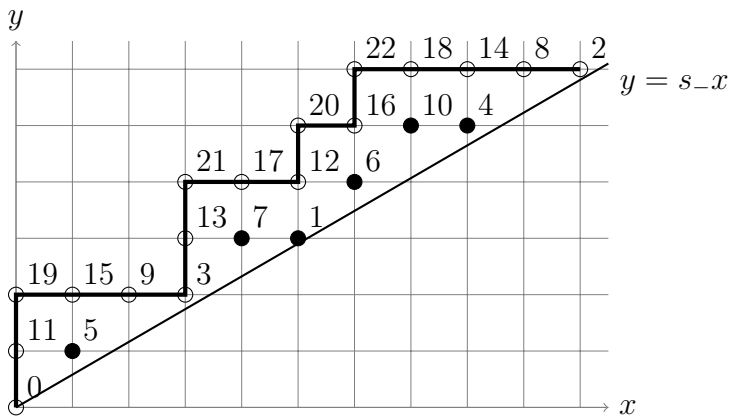
\begin{figure}
\begin{tikzpicture}[scale=.75]
\begin{scope}<+->;
  \draw[gray,step=1cm,very thin] (0,0) grid (10.5,6.5);

  \draw[gray,thin,->] (0, 0) -- (10.5, 0) node[black,right] {$x$};
  \draw[gray,thin,->] (0, 0) -- (0, 6.5) node[black,above] {$y$};
  
  \draw[black,thick,-] (0, 0) -- (10.5, 6.1) node[below right] {$y=s_- x$};
  \draw[black,ultra thick,-] (0, 0) -- (0, 2) -- (3, 2) -- (3, 4)
  -- (5, 4) -- (5,5) -- (6,5) -- (6,6) -- (10,6);
	\draw (0,0) node[draw,circle, inner sep=2pt] {} node[above right] {$0$};
	\draw (5,3) node[draw,circle, inner sep=2pt, fill] {} node[above right] {$1$};
	\draw (10,6) node[draw,circle, inner sep=2pt] {} node[above right] {$2$};
	\draw (3,2) node[draw,circle, inner sep=2pt] {} node[above right] {$3$};
	\draw (8,5) node[draw,circle, inner sep=2pt, fill] {} node[above right] {$4$};
	\draw (1,1) node[draw,circle, inner sep=2pt, fill] {} node[above right] {$5$};
	\draw (6,4) node[draw,circle, inner sep=2pt, fill] {} node[above right] {$6$};
	\draw (4,3) node[draw,circle, inner sep=2pt, fill] {} node[above right] {$7$};
	\draw (9,6) node[draw,circle, inner sep=2pt] {} node[above right] {$8$};
	\draw (2,2) node[draw,circle, inner sep=2pt] {} node[above right] {$9$};
	\draw (7,5) node[draw,circle, inner sep=2pt, fill] {} node[above right] {$10$};

	\draw (0,1) node[draw,circle, inner sep=2pt] {} node[above right] {$11$};
	\draw (5,4) node[draw,circle, inner sep=2pt] {} node[above right] {$12$};
	\draw (3,3) node[draw,circle, inner sep=2pt] {} node[above right] {$13$};
	\draw (8,6) node[draw,circle, inner sep=2pt] {} node[above right] {$14$};
	\draw (1,2) node[draw,circle, inner sep=2pt] {} node[above right] {$15$};
	\draw (6,5) node[draw,circle, inner sep=2pt] {} node[above right] {$16$};
	\draw (4,4) node[draw,circle, inner sep=2pt] {} node[above right] {$17$};
	\draw (7,6) node[draw,circle, inner sep=2pt] {} node[above right] {$18$};

	\draw (0,2) node[draw,circle, inner sep=2pt] {} node[above right] {$19$};
	\draw (5,5) node[draw,circle, inner sep=2pt] {} node[above right] {$20$};
	\draw (3,4) node[draw,circle, inner sep=2pt] {} node[above right] {$21$};
	\draw (6,6) node[draw,circle, inner sep=2pt] {} node[above right] {$22$};
\end{scope}
\end{tikzpicture}
\caption{A $(10,6)$-Dyck path.}
\label{fig:dp10x6}
\end{figure}

On Figure \ref{fig:dp10x6} we show an example of a $(10, 6)$-Dyck path. Encoding a North step by $1$ and an East step by $0$ the path on the picture is encoded by the sequence
\[
\pi=(1,1,0,0,0,1,1,0,0,1,0,1,0,0,0,0)\in\DD_{10, 6}.
\]
The \emph{reading order} is the order on the lattice points by their distance to the line $y=s_- x$. On the figure, next to each lattice point between the path and the diagonal, we put their position in the order. Recall that the \emph{area} statistic is defined as the number of cells between the path and the diagonal. This is clearly the same as the number of lattice points strictly between the path and the line $y=s_- x$. These points are marked on the Figure as $\bullet$:
\[
\area(\pi) = \#\{\bullet\}.
\]

A \emph{word parking function} is a map $w$ from the set of North steps to the set $\Z_{>0}$ such that for each two consecutive North steps $i$, $j$ the labels decrease: $w_i>w_j$. Denote the set of word parking functions by $\WP_\pi$. We will identify the North steps with the lattice points of their beginnings. For instance, on the Figure the North steps correspond to the points labeled $0, 3, 11, 12, 13, 16$.
We say that a North step attacks all the North steps which appear after it in the reading order, but before the lattice point directly above it. In our example $0$ attacks all points with labels between $0$ and $11$, $3$ attacks all points between $3$ and $13$, and so on. The \emph{temporary dinv} statistic of a word parking function is then defined in \cite{bergeron2015compositional} as
\[
\tdinv(\pi, w) = \#\{i,j:\; \text{$i$ attacks $j$ and $w_i>w_j$}\}
\]

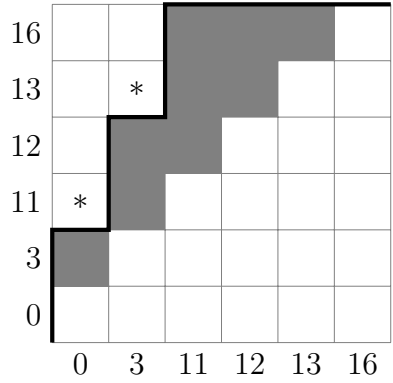
\begin{figure}
\begin{tikzpicture}[scale=.75]
\begin{scope}<+->;
  \draw[gray,step=1cm,very thin] (0,0) grid (6,6);

  \draw (0.5, 0) node[below] {$0$};
  \draw (1.5, 0) node[below] {$3$};
  \draw (2.5, 0) node[below] {$11$};
  \draw (3.5, 0) node[below] {$12$};
  \draw (4.5, 0) node[below] {$13$};
  \draw (5.5, 0) node[below] {$16$};
  \draw (0, 0.5) node[left] {$0$};
  \draw (0, 1.5) node[left] {$3$};
  \draw (0, 2.5) node[left] {$11$};
  \draw (0, 3.5) node[left] {$12$};
  \draw (0, 4.5) node[left] {$13$};
  \draw (0, 5.5) node[left] {$16$};
  \fill[gray] (0, 1) rectangle (1,2); 
  \fill[gray] (1, 2) rectangle (2,4); 
  \fill[gray] (2, 3) rectangle (3,6); 
  \fill[gray] (3, 4) rectangle (4,6); 
  \fill[gray] (4, 5) rectangle (5,6); 
  \draw (0.5, 2.5) node[black] {$*$};
  \draw (1.5, 4.5) node[black] {$*$};
  
  \draw[black,ultra thick,-] (0, 0) -- (0, 2) -- (1, 2) -- (1, 4)
  -- (2, 4) -- (2,6) -- (6,6);
\end{scope}
\end{tikzpicture}
\caption{The corresponding $\pi'$, $S_\pi$.}
\label{fig:dp10x6prime}
\end{figure}

Plotting the graph of the attack relation we obtain an $(n,n)$-Dyck path $\pi'$. Pairs of consecutive North steps form a subset of the set of corners of $\pi'$, which we denote $S_\pi$. Then $S_\pi$-admissible sequences are precisely the word parking functions, and we have (see Section \ref{sec:charfunc}) $\tdinv(\pi, w) = \inv(\pi', w)$, and
\[
\sum_{w\in\WP_\pi} q^{\tdinv(\pi, w)} x_w = \chi(\pi', S_\pi).
\]
In our example we obtain $\pi'$ as displayed on Figure \ref{fig:dp10x6prime}. The marked corners are denoted by $*$.

\begin{figure}
\begin{tikzpicture}[scale=.75]
  \draw[black,ultra thick,-] (0, 0) -- (1, 0);
  \draw[gray,very thin] (0, 3) rectangle (1,4);
  \draw[black,ultra thick,-] (5, 3) -- (5, 4);
  \draw[>=stealth,black,thick,<->] (0.5, 0) -- (0.5, 3);
  \draw[black] (0.5,1.5) node[right] {$l$};
  \draw[>=stealth,black,thick,<->] (1, 3.5) -- (5, 3.5);
  \draw[black] (3,3.5) node[above] {$a$};
  \draw (1, 0) node [above right] {$i$};
  \draw (5, 3) node [above right] {$j$};
  \draw[dashed] (1,0) -- (5,4);
  \draw[dashed] (0,0) -- (5,3);
\end{tikzpicture}
\caption{Arms and legs}
\label{fig:armleg}
\end{figure}
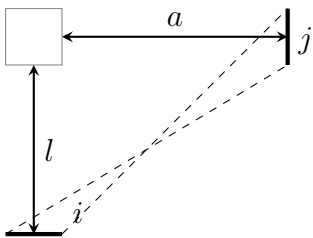

As explained in \cite{bergeron2015compositional}, the $\tdinv$ statistic needs to be modified. Let $\maxtdinv(\pi)$ be the maximal value $\tdinv$ can attain for $w\in\WP_\pi$. This is, clearly, equal to $\area(\pi')$. Let $\dinv(\pi)$ denote the \emph{dinv statistic} $\pi$. This is defined as the number of pairs $i$, $j$ where $i$ is an East step, $j$ is a North step, $i$ is to the left of $j$, and the following inequalities hold \footnote{It seems, in equation (2.12) of \cite{bergeron2015compositional}, the authors forgot that the first inequality should be made non-strict for the definition to be correct in the non-coprime case}
\[
\frac{a}{l+1} \leq \frac{m}{n} < \frac{a+1}{l},
\]
see Figure \ref{fig:armleg}. The condition is equivalent to $s_-$ being between the slopes of the dashed lines on the figure, that is the line connecting the beginnings of $i$ and $j$ and the line connecting the ends of $i$ and $j$. This is equivalent to existence of a line with slope $s_-$ which intersects both steps. We agree that the case $(a,l)=(0,0)$ always contributes to the dinv.  

Now the right hand side of Conjecture 3.3 of \cite{bergeron2015compositional} becomes
\begin{equation}\label{eq:rhs}
D_{n/m}^\alpha:=\sum_{\pi\in\DD_{n/m}^\alpha} t^{\area(\pi)} q^{\dinv(\pi)-\maxtdinv(\pi)} \chi(\pi', S_\pi),
\end{equation}
where for any composition $\alpha=(\alpha_1, \ldots, \alpha_r)$ with $\sum_{i=1}^r \alpha_i=\gcd(m,n)$ we denote by $\DD_{n/m}^\alpha$ the set of $(m,n)$-Dyck paths which touch the diagonal precisely at the points $(m_1 \sum_{j<i} \alpha_j, n_1 \sum_{j<i} \alpha_j)$ for $i=1,\ldots,r+1$.

\subsection{The sweep process}\label{sec:sweep}
For an $(m,n)$-Dyck path $\pi$ we concentrate on the corresponding summand of \eqref{eq:rhs}
\[
D_\pi := t^{\area(\pi)} q^{\dinv(\pi)-\maxtdinv(\pi)} \chi(\pi', S_\pi).
\]
Here we describe an algorithm to calculate $D_\pi$ by a process in which we move a line parallel to $s_-$ downwards, keeping track of what is happening on the line. This corresponds to following the path $\pi'$ as in Theorem \ref{thm:combinatorial}. The initial state corresponds to a line which is completely above $\pi$. We associate to the initial state the element $\varphi=1\in V_0$. As we are moving the line downwards we register an \emph{event} if the line passes through a lattice point weakly below the path. The final state will be the line $y=s_- x - \varepsilon$, so that all the points weakly between the path and the diagonal have been passed. In the final state $\varphi$ will be an element of $V_0$ equal to $D_\pi$. In some intermediate state $\varphi$ is an element of $V_k$, where $k$ equals the number of connected components of the intersection of the line with the figure bounded by the path and the diagonal. It remains to specify the operations we perform on $\varphi$ when we register an event:

\begin{tabular}{m{4cm} m{4cm} m{4cm}}
A) & 
\begin{tikzpicture}[scale=0.5]
\fill[gray] (0,-1) rectangle (1,0);
\draw[black,ultra thick](0,-1) -- (0,0) -- (1,0);
\draw (0,0) node[draw,circle, inner sep=2pt] {};
\draw[dashed](-1,-1)--(1,1);
\draw (-1.5, -1.5) node {};
\draw (1.5, 1.5) node {};
\end{tikzpicture}
&
apply $d_+$\\
B) & 
\begin{tikzpicture}[scale=0.5]
\fill[gray] (-1,-1) rectangle (0,0);
\fill[gray] (0,-1) rectangle (1,0);
\fill[gray] (0,0) rectangle (1,1);
\draw[black,ultra thick](-1,0) -- (0,0) -- (0,1);
\draw (0,0) node[draw,circle, inner sep=2pt] {};
\draw[dashed](-1,-1)--(1,1);
\draw (-1.5, -1.5) node {};
\draw (1.5, 1.5) node {};
\end{tikzpicture}
&
apply $d_-$\\
C) & 
\begin{tikzpicture}[scale=0.5]
\fill[gray] (0,-1) rectangle (1,0);
\fill[gray] (0,0) rectangle (1,1);
\draw[black,ultra thick](0,-1) -- (0,0) -- (0,1);
\draw (0,0) node[draw,circle, inner sep=2pt] {};
\draw[dashed](-1,-1)--(1,1);
\draw (-1.5, -1.5) node {};
\draw (1.5, 1.5) node {};
\end{tikzpicture}
&
apply $q^{-a} \frac{d_- d_+ - d_+ d_-}{q-1}$\\
D) & 
\begin{tikzpicture}[scale=0.5]
\fill[gray] (-1,-1) rectangle (0,0);
\fill[gray] (0,-1) rectangle (1,0);
\draw[black,ultra thick](-1,0) -- (0,0) -- (1,0);
\draw (0,0) node[draw,circle, inner sep=2pt] {};
\draw[dashed](-1,-1)--(1,1);
\draw (-1.5, -1.5) node {};
\draw (1.5, 1.5) node {};
\end{tikzpicture}
&
multiply by $q^a$ \\
E) & 
\begin{tikzpicture}[scale=0.5]
\fill[gray] (-1,-1) rectangle (1,1);
\draw (0,0) node[draw,circle, inner sep=2pt,fill] {};
\draw[dashed](-1,-1)--(1,1);
\draw (-1.5, -1.5) node {};
\draw (1.5, 1.5) node {};
\end{tikzpicture}
&
multiply by $t$ \\
\end{tabular}

In the case of events C) and D) we denote by $a$ the number of vertical steps the line crosses to the right of the event. Note that in the very end, when we cross the point $(0,0)$ we have to apply B). We have
\begin{thm}
The above algorithm computes $D_\pi$.
\end{thm}
\begin{proof}
First we show that applying the operations A), B), and C) without the factor $q^{-a}$ computes $\chi(\pi', S_\pi)$. At some point of time let $k$ be the number of connected components of the intersection of the line with the figure bounded by the path and the diagonal. This equals to the number of intersections of the line with a North step. Registering event A) means that we have a new North step which attacks all the $k$ North steps we already have. Thus $\pi'$ is extended by a horizontal step and we apply $d_+$. Registering event B) means that one North step does not intersect our line anymore, so all the subsequent North steps do not attack it. Thus $\pi'$ is extended by a vertical step and we apply $d_-$. Finally, event C) means that one North step disappears and another one appears attacking all the North steps intersecting the line. Because the new step is directly below the old one, we have to add a marked corner to $\pi'$, so we apply $\frac{d_- d_+ - d_+ d_-}{q-1}$.

It is clear that the effect of the rule E) on the final result is precisely the factor $t^{\area(\pi)}$.

Finally, we have to make sure that the power of $q$ produced by the rules C) and D) equals $\dinv(\pi) - \maxtdinv(\pi)$. For this we note that $\maxtdinv$ equals to the number of pairs of North steps which are in the attack relation. This is the same as the number of pairs $(i,j)$ of North steps such that there exists a line with slope $s_-$ which intersects both. After A) we have an extra North step and an extra East which change the $\dinv$ by $+k$. Also $\maxtdinv$ is changed by $+k$, so the net effect is $0$. After B) and E) there are no extra steps, so the numbers do not change. After C) we have an extra North step, which gives $+(k-a)$ to $\dinv$ (the number of East steps to the left intersecting the line), while giving $+k$ to $\maxtdinv$. So the net effect is $-a$. After D) we have an extra East step, which gives $+a$ to $\dinv$ (the number of North steps to the right intersecting the line), and no extra North steps. Thus the net change is $+a$.
\end{proof}

\subsection{The recursion}
Before a rigorous formulation we explain the main idea.
Having an algorithm to compute $D_\pi$, it is more or less straightforward to produce a recursive procedure to compute the sum of $D_\pi$ over all $\pi$. The idea comes from Dynamic Programming. Note that all the operators in the rules A) -- E) are linear. So if we have to apply same sequences of operators to different arguments and then add the results, we can instead add the arguments first, and then apply the operators. Hopefully, the same sums of arguments can be reused, so we won't need to recompute them. 
An intermediate state of our algorithm is encoded by the following object.

\begin{defn}\label{defn:colorings}
Let $h\in \R_{> 0}$ be such that the line $l_h:=\{x,y:y=s_- x + h\}$ does not pass through a lattice point in the rectangle $[0,m]\times[0,n]$. An \emph{admissible coloring} of $l_h$ is a finite family of non-overlapping intervals $[a_i, b_i]\subset l_h\cap [0,m]\times[0,n]$ such that each $a_i$ belongs to a vertical lattice line, and each $b_i$ belongs to a horizontal lattice line, i.e. $x(a_i)\in\Z$ and $y(b_i)\in\Z$ for all $i$.
\end{defn}

It is clear that admissible colorings of $l_h$ are precisely the sets that can appear as the intersection of $l_h$ with the figure bounded by some $(m,n)$-Dyck path and the diagonal. Denote the set of admissible colorings by $\CC_h$. For each $(m,n)$-Dyck path $\pi$ the sweep process produces a family of admissible colorings which begins with the empty coloring for $h=h_0=n+\varepsilon$ and ends with a certain coloring when $h=h_1=n-m s_- + \varepsilon$. The final coloring precisely corresponds to the touch composition of $\pi$. The sums over Dyck paths become sums over possible families of colorings.

For each $h$ and $c\in\CC_h$ denote by $D_{s, c}\in V_k$ the sum of the intermediate results\footnote{These elements can be interpreted as sums over appropriately defined ``partial parking functions'', see definition of a partial Dyck path and its character in \cite{carlsson2015proof} for details.} of the algorithm of Section \ref{sec:sweep} over all possible families of colorings which begin with the empty coloring and end at the coloring $c$. Here $k$ is the number of components of $c$. Now we reverse the time and find a recursion for $D_{s, c}$. So we start with the line $l=l_h$, and move it upwards until we hit a lattice point $P=(x_0, y_0)$ inside $c$ or on the boundary. Assume $h$ is such that $l_h$ passes through $P$, and denote $h_- = h-\varepsilon$, $h_+=h+\varepsilon$, so $c\in \CC_{h_-}$. If $P$ is on the boundary of $c$ at the moment $h$, this means $c$ comes from a unique $c'\in \CC_{h_+}$ by the rules A), C) or D). If $P$ is inside $c$, then $c$ can be obtained in $2$ ways: from $c'\in \CC_{h_+}$ by rule B) or from $c''\in \CC_{h_+}$ by rule E). 

Thus we obtain the following
\begin{thm}\label{thm:coloringrec}
For each $h\geq h_1$ and $c\in \CC_h$ we have an element $D_{s,c}\in V_k$, where $k$ is the number of components of $c$. These elements can be computed by the following recursion rules starting from $D_{s, \emptyset} = 1\in V_0$:

\begin{tabular}{m{2cm} m{5cm} m{6cm}}
A) & 
$D_{s,c} = d_+ D_{s,c'}$ &
\begin{tikzpicture}[scale=0.75]
\fill[gray] (0,-1) rectangle (1,0);
\draw[black,ultra thick](0,-1) -- (0,0) -- (1,0);
\draw (0,0) node[draw,circle, inner sep=2pt] {};
\draw[dashed](-0.5,-1)--(1, 0.5) node[above right] {$c$};
\draw[dashed](-1,-0.5)--(0.5, 1) node[above right] {$c'$};
\draw (-1.5, -1.5) node {};
\draw (1.5, 1.5) node {};
\end{tikzpicture}
\\
C) & 
$D_{s,c} = q^{-a} \frac{d_- d_+ - d_+ d_-}{q-1} D_{s,c'}$
&
\begin{tikzpicture}[scale=0.75]
\fill[gray] (0,-1) rectangle (1,0);
\fill[gray] (0,0) rectangle (1,1);
\draw[black,ultra thick](0,-1) -- (0,0) -- (0,1);
\draw (0,0) node[draw,circle, inner sep=2pt] {};
\draw[dashed](-0.5,-1)--(1, 0.5) node[above right] {$c$};
\draw[dashed](-1,-0.5)--(0.5, 1) node[above right] {$c'$};
\draw (-1.5, -1.5) node {};
\draw (1.5, 1.5) node {};
\end{tikzpicture}
\\
D) & 
$D_{s,c} = q^{a} D_{s,c'}$
&
\begin{tikzpicture}[scale=0.75]
\fill[gray] (-1,-1) rectangle (0,0);
\fill[gray] (0,-1) rectangle (1,0);
\draw[black,ultra thick](-1,0) -- (0,0) -- (1,0);
\draw (0,0) node[draw,circle, inner sep=2pt] {};
\draw[dashed](-0.5,-1)--(1, 0.5) node[above right] {$c$};
\draw[dashed](-1,-0.5)--(0.5, 1) node[above right] {$c'$};
\draw (-1.5, -1.5) node {};
\draw (1.5, 1.5) node {};
\end{tikzpicture}
\\
BE) & 
$D_{s,c} = d_- D_{s,c'} + t D_{s, c''}$
&
\begin{tikzpicture}[scale=0.75]
\fill[gray] (-1,-1) rectangle (0,0);
\fill[gray] (0,-1) rectangle (1,0);
\fill[gray] (0,0) rectangle (1,1);
\draw[black,ultra thick](-1,0) -- (0,0) -- (0,1);
\draw (0,0) node[draw,circle, inner sep=2pt] {};
\draw[dashed](-0.5,-1)--(1, 0.5) node[above right] {$c$};
\draw[dashed](-1,-0.5)--(0.5, 1) node[above right] {$c'$};
\draw (-1.5, -1.5) node {};
\draw (1.5, 1.5) node {};
\end{tikzpicture}
\begin{tikzpicture}[scale=0.75]
\fill[gray] (-1,-1) rectangle (1,1);
\draw (0,0) node[draw,circle, inner sep=2pt,fill] {};
\draw[dashed](-0.5,-1)--(1, 0.5) node[above right] {$c$};
\draw[dashed](-1,-0.5)--(0.5, 1) node[above right] {$c''$};
\draw (-1.5, -1.5) node {};
\draw (1.5, 1.5) node {};
\end{tikzpicture}
\\
\end{tabular}
\end{thm}

\begin{rem}\label{rem:compositioncase}
Note that for any composition $\alpha=(\alpha_1, \alpha_2,\ldots, \alpha_r)$ with $\sum_{i=1}^r \alpha_i = \gcd(m,n)$ we have
\begin{equation}\label{eq:compstocolorings}
D_{n/m}^\alpha = t^{\sum_{i=1}^r(\alpha_i-1)} d_-^{r} D_{n/m, c_\alpha},
\end{equation}
where $c_\alpha$ is the union of $r$ segments such the $i$-th segment begins at $x_i=m_1\sum_{j<i} \alpha_j$ and ends at $y_i=n_1 \sum_{j\leq i} \alpha_j$. This is because by passing from $h=h_1$ to $h=-\varepsilon$ we need to apply rule B) $r$ times. The extra power of $t$ comes from the $\sum_{i=1}^r(\alpha_i-1)$ points on the diagonal $y=s x$ that were not counted earlier.
\end{rem}

\section{Braids}
In order to proceed, we need to find a good expression for $D_{s, c}$ using the generators of $\bA_{q,t}$, and then prove that it satisfies the same recursion relations. It turns out, that the required expression is simply given by the braid obtained by wrapping $c$ around the punctured torus $(\R/\Z)^2\setminus\{(0,0)\}$. 

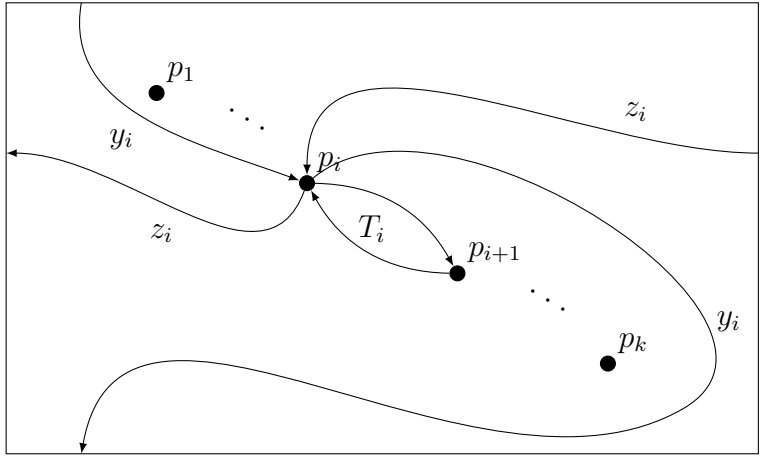
\begin{figure}
\begin{tikzpicture}[scale=2]
\draw (0,0) rectangle (5,3);
\draw (1.0, 3-0.6) node[draw,circle, inner sep=2pt,fill] {} node[above right] {$p_1$};
\draw (1.5-0.1, 3-0.9+0.06) node[above right] {$\cdot$};
\draw (1.5, 3-0.9) node[above right] {$\cdot$};
\draw (1.5+0.1, 3-0.9-0.06) node[above right] {$\cdot$};
\draw (2, 3-1.2) node[draw,circle, inner sep=2pt,fill] (1) {} node[above right] {$p_i$};
\draw (3, 1.2) node[draw,circle, inner sep=2pt,fill] (2) {} node[above right] {$p_{i+1}$};
\draw (3.5-0.1, 0.9+0.06) node[above right] {$\cdot$};
\draw (3.5, 0.9) node[above right] {$\cdot$};
\draw (3.5+0.1, 0.9-0.06) node[above right] {$\cdot$};
\draw (4,0.6) node[draw,circle, inner sep=2pt,fill] {} node[above right] {$p_k$};
\draw[>=latex,thin,->] (1) to [bend right=-30]  node[midway, below left] {$T_i$} (2); 
\draw[>=latex,thin,->] (2) to [bend right=-30] (1); 
\draw[>=latex,thin,<-] (1) to [out=90, in=180]  node[near end, above right] {$z_i$} (5, 3-1); 
\draw[>=latex,thin,<-] (0, 3-1) to [out=0, in=-110]  node[midway, below left] {$z_i$} (1); 
\draw[>=latex,thin,->] (1) to [out=40, in=30]  node[near end, right] {$y_i$} (4.5, 0.3)
to [out=-150, in=80] (0.5,0);
\draw[>=latex,thin,->] (0.5, 3) to [out=-100, in=160]  node[midway, below left] {$y_i$} (1); 
\end{tikzpicture}
\caption{The base point and the generators of $\B_k(\TT_0)$}
\label{fig:generators}
\end{figure}

We begin by giving an explicit presentation of the braid group of punctured torus. Denote by $\TT_0:=(\R/\Z)^2\setminus\{(0,0)\}$ the punctured torus. Its $k$-th braid group $\B_k(\TT_0)$ is the fundamental group of the space of $k$ distinct points on $\TT_0$. We fix a base point of this space by placing $k$ points along the diagonal with slope $-1$, see Figure \ref{fig:generators}. The torus is represented as a rectangle with the opposite sides glued together. The points $p_1,\ldots, p_k$ are situated along the diagonal. For each $i=1,\ldots, k-1$ the generator $T_i$ permutes $p_i$ and $p_{i+1}$ rotating them clockwise. For each $i=1,\ldots, k$ the generator $z_i$ moves $p_i$ along the vector $(-1,0)$ for one turn around the torus. The generator $y_i$ moves $p_i$ upwards and to the right around the other points, then downwards around the torus, then to the left around the other points to the original location.

After minor modifications, we have:
\begin{thm}[\cite{bellingeri2004presentations}\footnote{Strictly speaking, our results are independent of this Theorem, nevertheless, we find it important for understanding.}]\label{thm:braid_presentation}
For $k=1$ the group $\B_k(\TT_0)$ is the free group generated by $y_1$ and $z_1$. For $k>1$ it is the group generated by $T_i$ for $i=1,\ldots,k-1$ and $y_i$, $z_i$ for $i=1,\ldots,k$ subject to the relations
\[
T_i T_{i+1} T_i = T_{i+1} T_i T_{i+1}, \quad T_i T_j = T_j T_i \quad (|i-j|>1),
\]
\[
y_i T_j = T_j y_i, \quad  z_i T_j = T_j z_i \quad (i\notin \{j, j+1\}),
\]
\[
y_{i+1} = T_i^{-1} y_{i} T_i^{-1}, \quad z_{i+1} = T_i z_i T_i \quad (1\leq i\leq k-1),
\]
\[
y_i y_j = y_j y_i, \quad z_i z_j = z_j z_i \quad (1\leq i,j\leq k),
\]
\[
z_1 T_1 y_1 T_1^{-1} = T_1^{-1} y_1 T_1^{-1} z_1.
\]
\end{thm}

We will map $y_i$, $z_i$ to the corresponding generators of $\bA_{q,t}$, and our choice of the generators will make sure that they satisfy the various relations that we will need later. Here is a vague idea behind this choice: The operator $d_-$ must commute with $z_i$ and $y_i$, so geometrically it should move the last point $p_k$ along a straight line towards $(0,0)$, and kill it there. On the other hand, $d_+$ should create a new point at $(0,1)$ and move it to the position of $p_1$, simultaneously shifting all the other points to the right. These properties can be translated to commutation relations between $d_+$ and $y_i$, $z_i$, and our choices make sure that these relations indeed hold.

The more natural choice where $y_i$ would simply move $p_i$ directly downwards does not satisfy our relations. We will denote the corresponding commuting elements by $\tilde y_i$:
\[
\tilde y_i := T_{i\downto 1}\, T_{1\upto k}\, y_k\, T_{k\downto i}^*,
\]
in particular,
\[
\tilde y_1 = T_{1\upto k}\, y_k\, T_{k\downto 1}^*,
\]
\[
\tilde y_k = T_{k\downto 1}\, y_1\, T_{1\upto k}^*.
\]

The relations of Theorem \ref{thm:braid_presentation} can be translated to the following relations for $\tilde{y_i}$:
\begin{equation}\label{eq:ytilde1}
\tilde y_i  T_j = T_j \tilde y_i  \quad (i\notin \{j, j+1\}), \quad \tilde y_{i+1}  = T_i \tilde y_i T_i \quad(1\leq i\leq k-1),
\end{equation}
\begin{equation}\label{eq:ytilde2}
\tilde{y_i} \tilde{y_j} = \tilde{y_j} \tilde{y_i}  \quad (1\leq i,j\leq k), \quad \tilde y_1 T_1 z_1 = T_1 z_1 T_1 \tilde y_1 T_1.
\end{equation}
We will often use the following identities, which show various interactions between the ``trains'' of the form $T_{a\upto b}$, $T_{a\downto b}$. Recall that $T_{a\upto b}$ for $a>b$ stands for $T_{a\downto b}^*$, and similarly for the opposite direction. For any $a$, $b$, and $c\neq b$ we denote
\[
\sigma_{a,b}(c) = \begin{cases}
c+1 & \text{if $a\leq c< b$,}\\
c-1 & \text{if $a\geq c > b$,}\\
c & \text{otherwise.}
\end{cases}
\]
\begin{lem}\label{lem:trains} We have the following commutation rules:
\[
\text{gluing:}\quad T_{a\upto b}\, T_{b\upto c} = T_{a\upto c} \quad \text{for any $a,b,c$;}
\]
\[
\text{collision:}\quad T_{a\upto b}\, T_{c\downto d} = T_{c'\downto d'}\;
T_{a'\upto b'} \quad \text{if $b\neq c$,}
\]
where 
\[
b'=\sigma_{d,c}(b),\; c'=\sigma_{a,b}(c),\; a'=\sigma_{d,c'}(a),\;d'=\sigma_{a',b'}(d);
\]
\[
\text{overtaking:}\quad T_{a\upto b} T_{c\upto d} = T_{c\upto d} T_{a'\upto b'}\quad \text{if $a,b\in [d,c)$,}
\]
where 
\[
b'=\sigma_{d,c}(b),\; a'=\sigma_{d,c}(a);
\]
\[
\text{$T$--$z$ and $T$--$\tilde y$ rules:}\quad T_{a\downto b} z_b = z_a T_{a\upto b},\quad T_{a\downto b}\; \tilde y_b = \tilde y_a T_{a\upto b} \quad \text{for any $a,b$}.
\]
\end{lem}
\begin{proof}
A tedious case-by-case analysis.
\end{proof}

\begin{prop}
Denote by $\B_k^+(\TT_0)$ the monoid generated by $T_i$ and $T_i^{-1}$ for $1\leq i \leq k-1$, $y_1$ and $z_1$, subject to the relations of Theorem \ref{thm:braid_presentation}. We have a representation $\B_k^+(\TT_0)\to \bA_{q,t}[q^{\frac12}]$ given by sending $T_i$ to $q^{-\frac12} T_i$, $y_i$ to $-y_i$ and $z_i$ to $(qt)^{-1} z_i$.
\end{prop}
\begin{proof}
Only the last relation of Theorem \ref{thm:braid_presentation} is not immediate. This, in $\bA_{q,t}$, translates to 
\[
z_1 T_1 y_1 T_1^{-1} = q T_1^{-1} y_1 T_1^{-1} z_1,
\]
which, expressing $z_1$ in terms of $d_-$ and $d_+^*$ and using the definition of $y_2$ becomes
\[
(d_+^* d_- - d_- d_+^*)\, T_{k\downto 2}^*\, y_1 T_1^{-1} = y_2 (d_+^* d_- - d_- d_+^*)\, T_{k\downto 1}^*,
\]
which follows from the relation $d_+^* y_1 = y_2 d_+^*$, a part of Definition \ref{defn:correctinter}.
\end{proof}
\begin{rem}\label{rem:braid_rels}
It is not clear if the natural map $\B_k^+(\TT_0) \to \B_k(\TT_0)$ is an inclusion, in other words, if the relations of Theorem \ref{thm:braid_presentation} are enough to generate all relations of the submonoid of $\B_k(\TT_0)$ generated by $T_i$, $T_i^{-1}$, $y_i$ and $z_i$.
\end{rem}

\subsection{Special braids}
\begin{figure}
\begin{tikzpicture}[scale=1.5]
\draw[gray,thin,->] (0, 0) -- (3.5, 0) node[black,right] {$x$};
\draw[gray,thin,->] (0, 0) -- (0, 3.5) node[black,above] {$y$};
\draw[dashed] (0,3) -- (3,0);
\draw (0,0) rectangle (3,3);
\draw (0.5,2.5) node[draw,circle, inner sep=2pt,fill] {} node[above right] {\tiny $x=x_1$};
\draw (0.5,2.5) -- (0,1.5);
\draw (3,1.5) -- (2.5,0.5) node[draw,circle, inner sep=2pt,fill] {} node[right] {\tiny $\;x=\opnext(x_1)$};

\draw (2,1) node[draw,circle, inner sep=2pt,fill] {} node[above right] {\tiny $x=x_2$};
\draw (2,1) -- (1.5,0);
\draw (1.5,3) -- (1,2) node[draw,circle, inner sep=2pt,fill] {} node[right] {\tiny $\;x=\opnext(x_2)$};
\end{tikzpicture}
\caption{The operation $\opnext(x)$}
\label{fig:opnext}
\end{figure}
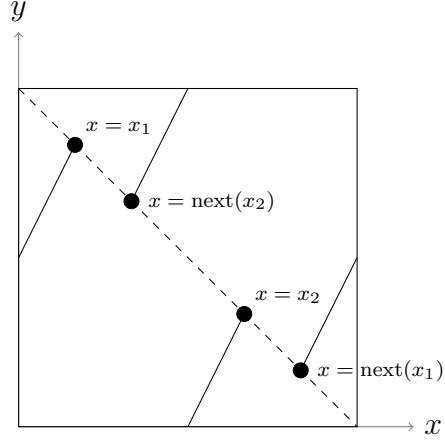

Fix $s\in\R_{>0}$. Let $t=1/(s+1)$. Let $I=(0,1)$ be the open unit interval. For $x\in I$ such that $x\neq t$ define 
\[
\opnext(x) := \begin{cases}
x-t & (x>t),\\
x+(1-t) & (x<t).
\end{cases}
\]
If we start at a point $(x,1-x)\in\TT_0$ and move downwards along the line with slope $s$, we will hit the line $y=1-x$ at the point $(\opnext(x), 1-\opnext(x))$. If $x<t$, we will have to cross the vertical wall. If $x>t$, we will cross the horizontal wall (see Figure \ref{fig:opnext}).

Now take $k$ distinct points $v=(v_1,v_2,\ldots, v_k)\in I^k$. For any $1\leq i\leq k$ such that $v_i\neq t$ and $\opnext(v_i)\neq v_j$ for all $j\neq i$ set 
\[
\opnext_i(v)_j = \begin{cases}
v_j & (j\neq i),\\
\opnext(v_i) & (j=i).
\end{cases}
\]
To the data of $v$, $i$ we also associate an element $b_i(v)\in\B_k^+(\TT_0)$ as follows:
\begin{equation}\label{eq:braidnext}
b_i(v) = \begin{cases}
T_{a'\downto a}\, z_a & (v_i<t),\\
T_{a'\upto a}^*\, \tilde y_a& (v_i>t),
\end{cases}
\end{equation}
where $a$ denotes the position of $v_i$ in the list of the elements of $v$ sorted in the ascending order, and similarly $a'$ is the position of $\opnext(v_i)$ in the sorted list of the elements of $\opnext_i(v)$.
Now to a sequence of indices $i_1,i_2,\ldots,i_l$ with $l\geq 0$ and a staring position $v$
we define inductively a position $\opnext_{i_1,i_2,\ldots, i_l}(v)$ and a braid $b_{i_1,i_2,\ldots, i_l}(v)\in \B_k^+(\TT_0)$ in such a way that
\[
\opnext_{\emptyset} (v) = v, \quad b_{\emptyset}(v) = \Id,
\]
\[
\opnext_{i_1,i_2,\ldots, i_l} (v) = \opnext_{i_1} \opnext_{i_2,\ldots, i_l} (v),\quad b_{i_1,i_2,\ldots, i_l} (v) = b_{i_1} (\opnext_{i_2,\ldots, i_l} (v)) b_{i_2,\ldots, i_l} (v),
\]
provided the points do not collide with each other and with the point $0$. A sequence $i_1, i_2, \ldots, i_l$ such that in $v$, $\opnext_{i_l} v$, $\opnext_{i_{l-1} i_l} v$, \dots no two points coincide and no point equals $0$ will be called \emph{admissible}

Geometrically, each point travels around the torus along the line of slope $s$ for some time. If the trajectories of the points do not intersect, it is evident that the resulting element in the braid group does not depend on the order of indices $i_1, i_2, \ldots, i_l$. Nevertheless, because of Remark \ref{rem:braid_rels} we give a direct combinatorial proof here.
\begin{lem}\label{lem:specialbraids}
Suppose the sequence $i_1, i_2, \ldots, i_l$ and all its permutations are admissible. Then the element $b_{i_1, i_2, \ldots, i_l}(v)\in\B_k^+(\TT_0)$ does not depend on the order of indices $i_1, i_2, \ldots, i_l$.
\end{lem}  
\begin{proof}
It is enough to permute only two indices, so we assume $l=2$. Suppose the positions of $i_1$, $i_2$ in the sorted list $v$ are $a_1$, $a_2$ respectively. Without loss of generality we assume $a_1<a_2$. Suppose after $\opnext_{i_1,i_2}$ the positions are $a_1'$, $a_2'$. We perform a case-by-case analysis as follows:
\begin{enumerate}
\item $(v_{i_1} < v_{i_2} < t) \;\Rightarrow\; a_1\leq a_1'<a_2',\; a_1<a_2\leq a_2'$,
\[
b_{i_1 i_2} (v) = T_{a_1'\downto a_1}\, z_{a_1} \,T_{a_2'\downto a_2}\, z_{a_2},
\]
\begin{enumerate}
\item $(v_{i_1} +1-t < v_{i_2}) \;\Rightarrow\; a_1'<a_2,$
\[
b_{i_2 i_1} (v) = T_{a_2'\downto a_2}\, z_{a_2}\, T_{a_1'\downto a_1}\, z_{a_1},
\]
\item $(v_{i_1} +1-t > v_{i_2}) \;\Rightarrow\; a_2\leq a_1',$
\[
b_{i_2 i_1} (v) = T_{a_2'\downto a_2-1}\; z_{a_2-1}\; T_{a_1'+1\downto a_1}\, z_{a_1},
\]
\end{enumerate}
\item $(t<v_{i_1} < v_{i_2}) \;\Rightarrow\; a_1'\leq a_1<a_2,\; a_1'<a_2'\leq a_2$,
\[
b_{i_2 i_1} (v) = T_{a_2'\upto a_2}^* \,\tilde y_{a_2}\, T_{a_1'\upto a_1}^* \, \tilde y_{a_1},
\]
\begin{enumerate}
\item $(v_{i_2} - t > v_{i_1}) \;\Rightarrow\; a_1<a_2',$
\[
b_{i_1 i_2} (v) = T_{a_1'\upto a_1}^* \, \tilde y_{a_1}\,T_{a_2'\upto a_2}^* \,\tilde y_{a_2},
\]
\item $(v_{i_2} -t < v_{i_1}) \;\Rightarrow\; a_2'\leq a_1,$
\[
b_{i_1 i_2} (v) = T_{a_1'\upto a_1+1}^*\; \tilde y_{a_1+1}\; T_{a_2'-1\upto a_2}^*\, \tilde y_{a_2},
\]
\end{enumerate}
\item $(v_{i_1}<t<v_{i_2}) \;\Rightarrow\; a_2' < a_1',\; a_1<a_1',\; a_2'< a_2$,
\begin{enumerate}
\item $(v_{i_2}<v_{i_1}+1-t) \;\Rightarrow\; a_2 \leq a_1',$
\[
b_{i_2 i_1} (v) = T_{a_2'\upto a_2-1}^*\; \tilde y_{a_2-1}\; T_{a_1'\downto a_1}\, z_{a_1},
\]
\item $(v_{i_2}>v_{i_1}+1-t) \;\Rightarrow\; a_1' \leq a_2,$
\[
b_{i_2 i_1} (v) = T_{a_2'\upto a_2}^*\, \tilde y_{a_2}\; T_{a_1'-1\downto a_1}\, z_{a_1},
\]
\end{enumerate}
\begin{enumerate}
\item[(a')] $(v_{i_2}-t>v_{i_1}) \;\Rightarrow\; a_1\leq a_2',$
\[
b_{i_1 i_2} (v) = T_{a_1'\downto a_1}\, z_{a_1}\; T_{a_2'+1\upto a_2}^*\, \tilde y_{a_2},
\]
\item[(b')] $(v_{i_2}-t<v_{i_1}) \;\Rightarrow\; a_2'\leq a_1,$
\[
b_{i_1 i_2} (v) = T_{a_1'\downto a_1+1}\; z_{a_1+1}\; T_{a_2'\upto a_2}^*\, \tilde y_{a_2},
\]
\end{enumerate}
\end{enumerate}

Note that in the case (iii) it is impossible to have $a_2'>a_1'$ because that would mean $v_{i_1}+1-t<v_{i_2}-t$, so $v_{i_2}>v_{i_1}+1$, which is impossible.
The cases (i) and (ii) are similar. For instance, in (i) the index $a_1$ is the smallest, while $a_2'$ is the largest. So we can move $z_{a_1}$ to the right, use $T_{a_2'\downto a_2}\, z_{a_2}=z_{a_2'}\, T_{a_2'\downto a_2}^*$,
and then move $z_{a_2'}$ to the left on both sides. Then the statement follows from the identities of Lemma \ref{lem:trains}.

In the case (iii) we have in both cases (a), (b)
\[
b_{i_2 i_1}(v) = T_{a_2'\downto 1}\, T_{a_1'\downto 2}\;  \tilde y_1 T_1 z_1 \;T_{2\upto a_2}\, T_{1\upto a_1}.
\]
This transforms to 
\[
T_{a_2'\downto 1}\, T_{a_1'\downto 1}\;  z_1 T_1 \tilde y_1 \;T_{1\upto a_2}\, T_{1\upto a_1} = T_{a_1'\downto 1}\, T_{a_2'+1\downto 2}\;  z_1 T_1 \tilde y_1 \;T_{2\upto a_1+1}\, T_{1\upto a_2},
\]
which equals $b_{i_1 i_2}(v)$ in both cases (a'), (b').
\end{proof}

By the Lemma, the element $b_{i_1, i_2, \ldots, i_l}(v)$ depends only on the number of times we use each index. This can be encoded in a composition and we introduce
\begin{defn}
Let $v=(v_1, v_2,\ldots, v_k)$ be a tuple of distinct points in $(0,1)$. Let $s\in\R_{>0}$. Let $\alpha=(\alpha_1, \alpha_2, \ldots, \alpha_k)$ be a composition\footnote{This composition should not be confused with the one specifying the touch points of a Dyck path} such that the tuple $1^{\alpha_1-1} 2^{\alpha_2-1},\ldots,k^{\alpha_k-1}$ is admissible in the sense that all the points $\opnext^j(v_i)$ for $i=1,\ldots, k$, $j=0,\ldots, \alpha_i-1$ are distinct and different from $0$. The corresponding \emph{special braid} is defined as 
\[
B_{s,v,\alpha} := b_{1^{\alpha_1-1}, 2^{\alpha_2-1},\ldots,k^{\alpha_k-1}}(v) \in \B_k^+(\TT_0).
\]
\end{defn}

\subsection{Creation operators}
These operators control how the braid changes when we add extra points  which do not move. This affects the indices $a$, $a'$ in the rules \eqref{eq:braidnext}. 

\begin{prop}\label{prop:phiplus}
Let $\varphi_+:\B_{k}^+(\TT_0)\to \B_{k+1}^+(\TT_0)$ be the homomorphism which sends $T_i$ to $T_{i+1}$, $z_i$ to $z_{i+1}$ and $\tilde y_i$ to $\tilde y_{i+1}$ for all $i$. For a special braid $B\in \B_{k}^+(\TT_0)$ let $B' \in \B_{k+1}^+(\TT_0)$ be the braid obtained by adding a point to the left of all the points of $B$ and their future positions. Then we have
\[
B' = \varphi_+(B).
\]
\end{prop}
\begin{proof}
 Adding an extra fixed point to the left of all the other points and their future positions increases $a$ and $a'$ by $1$ in each application of the rules \eqref{eq:braidnext} during construction of $B$.
\end{proof} 

\begin{prop}\label{prop:phiminus}
Let $\varphi_-:\B_{k}^+(\TT_0)\to \B_{k+1}^+(\TT_0)$ be the homomorphism that sends $T_i$ to $T_{i}$, $z_i$ to $z_{i}$ and $y_i$ to $y_i$ for all $i$. For a special braid $B\in \B_{k}^+(\TT_0)$ let $B' \in \B_{k+1}^+(\TT_0)$ be the braid obtained by adding a point $X$ at $(t, 1-t)$. Suppose the position of $X$ in the sorted list of initial positions of $B'$ is $i$, and in the final list is $i'$. Then we have
\begin{equation}\label{eq:phiminus}
T_{k+1\downto i'}^*\; B' = \varphi_-(B)\; T_{k+1\downto i}^*.
\end{equation}
\end{prop}
\begin{proof}
It is enough to verify the statement for the elementary moves \eqref{eq:braidnext}. Suppose first, that
\[
B = T_{a'\downto a} \; z_a.
\]
We necessarily have $a<i$. So the right hand side of \eqref{eq:phiminus} becomes
\[
T_{a'\downto a}\; z_a T_{k+1\downto i}^*
\]
\[
=\begin{cases}
T_{k+1\downto i}^*\;T_{a'\downto a}\; z_a
& \text{if $i>a'$,} \\
T_{k+1\downto i-1}^*\;T_{a'+1\downto a}\; z_a
& \text{if $i\leq a'$.} \\
\end{cases}
\]
If $i'=i$, then $a'<i$ and to obtain $B'$ we don't need any changes:
\[
B' = T_{a'\downto a}\; z_a.
\]
If $i'=i-1$, them $a'\geq i$, and to obtain $B'$ we need to increase $a'$ by $1$:
\[
B' = T_{a'+1\downto a}\; z_a.
\]

Now consider the second case in \eqref{eq:braidnext}. We have 
\[
B = T_{a'\upto a}^*\; \tilde y_a.
\]
We have $i\leq a$ and we don't know if $i'=i$ or $i'=i+1$. We first inspect the right hand side of the statement:
\[
\varphi_-(B)\; T_{k+1\downto i}^*
= T_{a'\upto k+1}^* \; \tilde y_{k+1}\; T_k\; T_{k\downto a}^*\; T_{k+1\downto i}^*,
\]
where to compute $\varphi_-(\tilde y_a)$ we used $\varphi_-(\tilde y_k) = T_k^{-1} \tilde y_{k+1} T_k$ and the expression of $\tilde y_a$ in terms of $\tilde y_k$. We further simplify it to
\[
T_{a'\upto k+1}^*\; \tilde y_{k+1} \;T_{k\downto i}^*\; T_{k+1\downto a+1}^* = 
T_{a'\upto k+1}^* \;T_{k\downto i}^*\; T_{k+1\downto a+1} \; \tilde y_{a+1} 
\]
\[
=T_{a'\upto k}^* \;T_{k\downto a} \; T_{k+1\downto i}^*\; \tilde y_{a+1} 
=T_{a'\upto a}^* \; T_{k+1\downto i}^*\; \tilde y_{a+1}
\]
\[
=\begin{cases}
T_{k+1\downto i}^*\; T_{a'+1\upto a+1}^*\; \tilde y_{a+1} & (i\leq a'),\\
T_{k+1\downto i+1}^*\; T_{a'\upto a+1}^*\; \tilde y_{a+1} & (i\geq a').
\end{cases}
\]
If $i'=i$, then $i\leq a'$, and to obtain $B'$ we need to increase both $a$ and $a'$:
\[
B'= T_{a'+1\upto a+1}^*\; \tilde y_{a+1}.
\]
On the other hand, if $i'=i+1$, then we have $i\geq a'$, and to obtain $B'$ we increase only $a$:
\[
B'= T_{a'\upto a+1}^*\; \tilde y_{a+1}.
\]
So we see that these cases precisely match the two cases above.
\end{proof}

Finally, we have a result for a point at position $(1-t, t)$.
\begin{prop}\label{prop:phiplus2}
Let $\varphi_+^*:\B_{k}^+(\TT_0)\to \B_{k+1}^+(\TT_0)$ be the homomorphism that sends $T_i$ to $T_{i+1}$, $z_i$ to $z_{i+1}$ and $y_i$ to $y_{i+1}$ for all $i$. For a special braid $B\in \B_{k}^+(\TT_0)$ let $B' \in \B_{k+1}^+(\TT_0)$ be the braid obtained by adding a point $X$ at $(1-t, t)$. Suppose the position of $X$ in the sorted list of initial positions of $B'$ is $i$, and in the final list is $i'$. Then we have
\begin{equation}\label{eq:phiplus}
B'\; T_{i\downto 1} = T_{i'\downto 1} \; \varphi_+^*(B).
\end{equation}
\end{prop}
\begin{proof}
The proof is analogous to the proof of Proposition \ref{prop:phiminus}. In the first case we have $a'\geq i'$, and $i'=i$ or $i'=i-1$. In any case, the right hand side of \eqref{eq:phiplus} is
\[
T_{i'\downto i} \; T_{a'+1\downto a+1}\; z_{a+1}
=\begin{cases}
T_{a'+1\downto a+1}\; z_{a+1}\; T_{i'\downto 1} & \text{if  $a\geq i'$,} \\
T_{a'+1\downto a}\; z_{a}\; T_{i'+1\downto 1} & \text{if  $a\leq i'$,} \\
\end{cases}
\]
If $i'=i$ then also $a \geq i$, and to obtain $B'$ we need to increase both $a$ and $a'$, so we have
\[
B' = T_{a'+1\downto a+1}\; z_{a+1}.
\]
If $i'=i-1$ then $a\leq i'$, to obtain $B'$ we only increase $a'$:
\[
B' = T_{a'+1\downto a}\; z_{a}.
\]

In the second case we have $a'<i'$, and $i'=i$ or $i'=i+1$. We have $\varphi_+^*(\tilde y_1) = T_1^{-1} \tilde y_1 T_1$, so the right hand side of \eqref{eq:phiplus} can be written as
\[
T_{i'\downto 1} \; T_{a'+1\upto a+1}^*\; T_{a+1\downto 2}\; T_1^{-1} \tilde y_1\; T_{1\upto a+1} = T_{i'\downto 1} \; T_{a'+1\downto 2}\; T_1^{-1} \tilde y_1\; T_{1\upto a+1}
\]
\[
= T_{a'\downto 1}\; T_{i'\downto 2}\; \tilde y_1 \; T_{1\upto a+1}
= T_{a'\downto 1}\; \tilde y_1\; T_{i'\downto 2} \; T_{1\upto a+1}
\]
\[
=\begin{cases}
T_{a'\downto a}\;\tilde y_a\; T_{i'\downto 1} & \text{if $a\leq i'-1$}\\
T_{a'\downto a+1}\;\tilde y_{a+1}\; T_{i'-1\downto 1} & \text{if $a\geq i'-1$}\\
\end{cases}
\]
On the other hand, if $i'=i$ then $a<i'$ and we don't need to change $a$, $a'$ for $B'$:
\[
B' = T_{a'\upto a}^*\; \tilde y_a.
\]
If $i'=i+1$ then $a\geq i$ and we need to increase $a$ by $1$:
\[
B' = T_{a'\upto a+1}^*\; \tilde y_{a+1}.
\]
\end{proof}

\subsection{Braid associated to a coloring}\label{sec:braidcoloring}
The colorings from Definition \ref{defn:colorings} can be viewed as special braids with slope $s_-$ as follows. For $h\in\R_{>0}$ let $c$ be an admissible coloring of the line $l_h$ with intervals $[a_i, b_i]$. For each $i$ we project $[a_i, b_i]$ to the torus $\TT$. Let $v_i$ be the first intersection point of $[a_i, b_i]$ with the diagonal $y=1-x$ of $\TT$ when followed from $b_i$ to $a_i$. Let $\alpha_i$ be the total number of times $[a_i, b_i]$ intersects the diagonal. We define the braid associated to $c$ as $B_{s_-, v, \alpha}$. The condition of admissibility for colorings implies that in the initial position all the points are to the left of $1-t$, and in the final position to the left of $t$ (recall that $t=\frac{1}{s+1}$).

Next we analyse evolution of the braid when $c$ changes as in Theorem \ref{thm:coloringrec}. For colorings denoted as $c$, $c'$, $c''$ we denote the corresponding braids by $B$, $B'$, $B''$. They act on $V_*$ via the algebra $\bA_{q,t}$, and we are interested in their actions on the vectors of the form $d_+^k(1)$.

\subsubsection{Transformation A)}
This clearly corresponds to creation of one extra strand of length $0$ to the left of all the other strands. By Proposition \ref{prop:phiplus}
\[
\varphi_+(B') = B.
\]
We now compute the action of $\varphi_+$ on $y_1$:
\[
\varphi_+(y_1) = \varphi_+(T_{1\upto k}^*\; \tilde y_k\; T_{k\downto 1}) = T_1 y_1 T_1^{-1}.
\]
Thus we see that the following diagram is commutative:
\[
\begin{tikzcd}
V_k \arrow{r}{d_+} \arrow{d}{B'} & V_{k+1} \arrow{d}{\varphi_+(B')=B}\\
V_k \arrow{r}{d_+} & V_{k+1}
\end{tikzcd}
\]
Hence we have
\begin{equation}\label{eq:braidruleA}
B\; d_+^{k+1}(1) = d_+ B'\; d_+^k(1).
\end{equation}

\subsubsection{Transformation C)}
Here the rightmost strand of $B'$ makes one extra turn crossing the horizontal wall and ending up leftmost. Thus to obtain $B$ from $B'$ we need to compose $B'$ with the second part of \eqref{eq:braidnext} with $a=k$, $a'=1$. This is given by 
\[
T_{1\upto k}^*\; \tilde y_k = y_1 \;T_{1\upto k}^*,
\]
which is mapped to $q^{\frac{1-k}2}\frac{d_- d_+ - d_+ d_-}{q-1}$ in $\bA_{q,t}$, so we have a commutative diagram with $C:=\frac{d_- d_+ - d_+ d_-}{q-1}$
\[
\begin{tikzcd}
V_k \arrow{rr}{\Id} \arrow{d}{B'} & & V_{k} \arrow{d}{B}\\
V_k \arrow{rr}{q^{\frac{1-k}2}C} & & V_{k}
\end{tikzcd}
\]
This implies
\begin{equation}\label{eq:braidruleC}
B \;d_+^k(1) = q^{\frac{1-k}2}\; C\; B'\; d_+^k(1).
\end{equation}

\subsubsection{Transformation D)}
Here the strand with the rightmost initial position changes its initial position to the leftmost and acquires and extra turn crossing the vertical wall in the beginning of its existence. Thus we need to pre-compose $B'$ with the first part of \eqref{eq:braidnext} with $a=1$, $a'=k$, which is given by 
\[
T_{k\downto 1}\; z_1 = z_k\; T_{k\downto 1}^*.
\]
Thus we have the following expression for $B\; d_+^k(1)$:
\[
(qt)^{-1} \; q^{\frac{k-1}2} \;B'\; z_k\; T_{k\downto 1}^*\; d_+^k(1) = (qt)^{-1} \; q^{\frac{k-1}2} \;B'\; d_+^{k-1} z_1 d_+ (1)
\]
\[
= - q^{\frac{k-1}2} \;B'\; d_+^{k-1} y_1 d_+^* (1) = q^{\frac{k-1}2} \;B'\; d_+^{k} (1).
\]
This shows:
\begin{equation}\label{eq:braidruleD}
B\; d_+^k(1) = q^{\frac{k-1}2} \;B'\; d_+^{k} (1).
\end{equation}

\subsubsection{Transformations B) and E)}
\begin{figure}
\def\eps{0.04}
\def\s{1.5}
\def\t{(1/(1+\s))}
\begin{tikzpicture}[scale=3]
\draw (0.5,-0.2) node {$B$};
\draw[gray,thin,->] (0, 0) -- (1.1, 0) node[black,right] {$x$};
\draw[gray,thin,->] (0, 0) -- (0, 1.1) node[black,above] {$y$};
\draw (0,0) rectangle (1,1);
\draw ({\t+\eps+(\t+\eps) / \s}, 1) --  ({\t+\eps}, {1-\t-\eps}) node[midway, right]{$B_1$} node[draw,circle, inner sep=1pt,fill] {} node[right] {$X_1$};
\draw ({\t+\eps}, {1-\t-\eps}) -- ({\t+\eps - (1-\t-\eps) / \s}, 0);
\draw ({\t+\eps - (1-\t-\eps) / \s}, 1) -- (0, {1-(\s*(\t+\eps) - (1-\t-\eps))});
\draw (1, {1-(\s*(\t+\eps) - (1-\t-\eps))}) -- ({1-\t+\eps}, {\t-\eps})  node[draw,circle, inner sep=1pt,fill] {} node[right] {$X_2$};
\draw ({1-\t+\eps}, {\t-\eps}) -- ({1-\t+\eps-(\t-\eps) / \s}, 0)  node[midway, right]{$B_2$};
\end{tikzpicture}
\begin{tikzpicture}[scale=3]
\draw (0.5,-0.2) node {$B'$};
\draw[gray,thin,->] (0, 0) -- (1.1, 0) node[black,right] {$x$};
\draw[gray,thin,->] (0, 0) -- (0, 1.1) node[black,above] {$y$};
\draw (0,0) rectangle (1,1);
\draw ({\t-\eps+(\t-\eps) / \s}, 1) --  ({\t-\eps}, {1-\t+\eps}) node[midway, right]{$B_1$} node[draw,circle, inner sep=1pt,fill] {} node[right] {$X_1$};
\draw ({1-\t-\eps}, {\t+\eps}) node[draw,circle, inner sep=1pt,fill] {} node[right] {$X_2$};
\draw ({1-\t-\eps}, {\t+\eps}) -- ({1-\t-\eps-(\t+\eps) / \s}, 0)  node[midway, right]{$B_2$};
\end{tikzpicture}
\begin{tikzpicture}[scale=3]
\def\eps{0.03};
\draw (0.5,-0.2) node {$B''$};
\draw[gray,thin,->] (0, 0) -- (1.1, 0) node[black,right] {$x$};
\draw[gray,thin,->] (0, 0) -- (0, 1.1) node[black,above] {$y$};
\draw (0,0) rectangle (1,1);
\draw ({\t-\eps+(\t-\eps) / \s}, 1) --  ({\t-\eps}, {1-\t+\eps}) node[midway, right]{$B_1$} node[draw,circle, inner sep=1pt,fill] {} node[right] {$X_1$};
\draw ({\t-\eps}, {1-\t+\eps}) -- (0, {1-\t+\eps - \s*(\t-\eps)});
\draw (1, {1-\t+\eps - \s*(\t-\eps)}) -- ({1-((1-\t+\eps) / \s - (\t-\eps))}, 0);
\draw ({1-((1-\t+\eps) / \s - (\t-\eps))}, 1) -- ({1-\t-\eps}, {\t+\eps})  node[draw,circle, inner sep=1pt,fill] {} node[right] {$X_2$};
\draw ({1-\t-\eps}, {\t+\eps}) -- ({1-\t-\eps-(\t+\eps) / \s}, 0)  node[midway, right]{$B_2$};
\end{tikzpicture}
\caption{The transformations B) and E)}
\label{fig:braidsBE}
\end{figure}
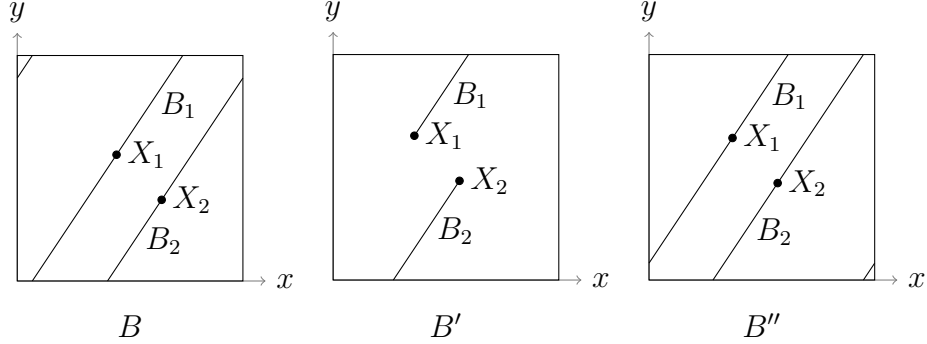

Next we obtain relations between the three braids $B, B'' \in \B_k^+(\TT_0)$ and $B'\in \B_{k+1}^+(\TT_0)$ from part BE of Theorem \ref{thm:coloringrec}. This case is more complicated than the previous ones as it involves $3$ braids. See Figure \ref{fig:braidsBE}.

We obtain $B$ as follows. Denote the component that passes next to the point $(0,0)$ by $i$. First we follow some of the components of $B$ different from $i$. Then we begin to follow $i$ and stop just before passing next to $(0,0)$. Denote the resulting braid by $B_1$ and the positions of the points by $v^{(1)}$. Denote the coordinate of $i$ at this moment by $X_1$ (see Figure \ref{fig:braidsBE}). To complete $B_1$ to $B$ we need to continue to cross the horizontal wall, then cross the vertical wall at the leftmost point, and then continue until the end following the remaining steps, which we denote by $B_2$. Denote by $v^{(2)}$ the positions of the points before performing the steps of $B_2$ and by $X_2$ the position of $i$. Suppose the positions of $i$ in the sorted $v^{(1)}$, $v^{(2)}$ are $a_1$, $a_2$ resp. Then we have
\[
B = B_2\; T_{a_2\downto 1}\; z_1 \;T_{1\upto a_1}^*\; \tilde y_{a_1} \;B_1
= B_2\; T_{a_2\downto 1}\; z_1 \tilde y_1 \;T_{1\upto a_1} \;B_1.
\]

The relative positions of points in $B''$ are the same as in $B$ at all steps except the two steps on Figure \ref{fig:braidsBE}. We have
\[
B'' = B_2\; T_{a_2\upto k}^* \;\tilde y_k\; T_{k\downto a_1}\; z_{a_1} \; B_1
= B_2\; T_{a_2\downto 1}\; y_1 z_1 \; T_{1\upto a_1} \;B_1
\]

Finally we consider the braid $B'$. First we perform the steps of $B_1$. However, the positions of points $a$, $a'$ in the rules \eqref{eq:braidnext} will be different because of presence of the point $X_2$. We apply Proposition \ref{prop:phiplus2} and obtain, noting that in the initial position all the points are to the left of $X_2$, that the first part of $B'$ equals:
\[
T_{a_2'\downto 1} \;\varphi_+^*(B_1)\; T_{1\upto k+1}^*,
\]
where $a_2'=a_2$ if $X_2<X_1$ and $a_2'=a_2+1$ if $X_2>X_1$.
After that, performing the steps of $B_2$ the indices $a$ and $a'$ will be affected by presence of $X_1$. Now we apply Proposition \ref{prop:phiminus} and obtain, noting that in the final position all the points are to the left of $X_1$, that the second part of $B'$ equals:
\[
\varphi_-(B_2) T_{k+1\downto a_1'}^*,
\]
where $a_1'=a_1$ if $X_1<X_2$ and $a_1'=a_1+1$ if $X_1>X_2$. It turns out that the following composition does not depend on whether $X_1>X_2$:
\[
T_{k+1\downto a_1'}^*\; T_{a_2'\downto 1}
=
\begin{cases}
T_{a_2'\downto 1} \;
T_{k+1\downto a_1'}^* & \text{if $a_1'>a_2'$,}\\
T_{a_2'-1\downto 1} \;
T_{k+1\downto a_1'+1}^* & \text{if $a_1'<a_2'$,}
\end{cases}
\;=T_{a_2\downto 1} \;
T_{k+1\downto a_1+1}^*.
\]
Therefore we have
\[
B' = \varphi_-(B_2) T_{a_2\downto 1} \;
T_{k+1\downto a_1+1}^* \;\varphi_+^*(B_1)\; T_{1\upto k+1}^*
\]
\[
= \varphi_-(\tilde B_2) T_{k+1\downto 2}^*\;\varphi_+^*(\tilde B_1) T_{1\upto k+1}^*,
\]
where
\[
\tilde B_1 = T_{1\upto a_1} \;B_1,\quad 
\tilde B_2 = B_2\; T_{a_2\downto 1}.
\]

Next we compute the evaluations at $d_+^k(1)$.
\[
B\, d_+^k (1) = -(qt)^{-1} \tilde B_2\, z_1 \tilde y_1 \, \tilde B_1\, d_+^k(1),
\]
\[
B''\, d_+^k (1) = -(qt)^{-1} \tilde B_2\, y_1 z_1 \, \tilde B_1\, d_+^k(1).
\]
For $B'$ using $T_i d_+^{k+1}(1) = d_+^{k+1}(1)$ and $d_+^k (1) = - y_1 d_+^* d_+^k (1)$ we obtain
\[
B' d_+^{k+1}(1) = - q^{\frac{2k-1}2} \varphi_-(\tilde B_2)\; T_{k+1\downto 2}^*\; \varphi_+^*(\tilde B_1) \,y_1\, d_+^* d_+^{k}(1).
\]
We have commutative diagrams, valid for any $B\in \B_k^+(\TT_0)$:
\[
\begin{tikzcd}
V_k \arrow{r}{y_1 d_+^*} \arrow{d}{B} & V_{k+1} \arrow{d}{\varphi_+^*(B)}\\
V_k \arrow{r}{y_1 d_+^*} & V_{k+1}
\end{tikzcd}
\begin{tikzcd}
V_k \arrow{r}{d_-} \arrow{d}{\varphi_-(B)} & V_{k+1} \arrow{d}{B}\\
V_k \arrow{r}{d_-} & V_{k+1}
\end{tikzcd}
\]
Thus we have
\[
d_- B' d_+^{k+1}(1) = - q^{\frac{2k-1}2} \tilde B_2\, d_- y_1 d_+^*\; T_{k\downto 1}^* \; \tilde B_1 d_+^k(1).
\]
The following identity in $\bA_{q,t}$ can be deduced from Definition \ref{defn:correctinter}:
\[
z_1 \tilde y_1 \;=\; t y_1 z_1 \;+\; t q^k d_- y_1 d_+^*\; T_{k\downto 1}^*.
\]
Thus we have
\begin{equation}\label{eq:braidruleBE}
B d_+^k(1) = t B'' d_+^k (1) + q^{-\frac12} d_- B' d_+^{k+1}(1).
\end{equation}

\subsection{The main result}
The expressions we have in \eqref{eq:braidruleA}--\eqref{eq:braidruleBE} differ slightly from the corresponding recursion relations in Theorem \ref{thm:coloringrec}. To obtain a precise match we need an extra factor which is a power of $q$. For an admissible coloring $c$ with $k$ components label the parts from left to right the by numbers from $1$ to $k$. Let $\sigma_{\mathrm{initial}}(c)$ be the permutation obtained by reading the labels in the initial position of the braid, and $\inv_{\mathrm{initial}}(c)$ be the number of inversions of $\sigma_{\mathrm{initial}}(c)$. Similarly define $\sigma_{\mathrm{final}}(c)$, $\inv_{\mathrm{final}}(c)$ for the final position.
\begin{thm}\label{thm:main}
For each slope $s$ and an admissible coloring $c$ with $k$ parts let $B_{s,c}\in\B_k^+(\TT_0)$ denote the corresponding special braid constructed in Section \ref{sec:braidcoloring}. Then the invariant $D_{s,c}\in V_k$ from Theorem \ref{thm:coloringrec} is related to $B_{s,c}$ as follows:
\[
D_{s,c} = q^{\frac{\inv_\mathrm{final}(c) - \inv_\mathrm{initial}(c)}2}\; B_{s,c}\; d_+^k (1).
\]
\end{thm} 
\begin{proof}
After we established \eqref{eq:braidruleA}--\eqref{eq:braidruleBE} it only remains to show that the extra factor $q^{\frac{\inv_\mathrm{final}(c) - \inv_\mathrm{initial}(c)}2}$ fixes the powers of $q$ correctly. So we consider the cases A)-E) one by one:
\begin{itemize}
\item[A)] Let $a$ be the number of parts to the left of the new part of $c$. We have 
\[
\inv_\mathrm{initial}(c) = \inv_\mathrm{initial}(c') + a,\quad \inv_\mathrm{final}(c) = \inv_\mathrm{final}(c') + a,
\]
as the new part has label $a+1$ and is the first in both $\sigma_{\mathrm{initial}}(c)$, $\sigma_{\mathrm{final}}(c)$.
\item[C)] Let $a$ be the number of parts to the right of the part of $c$ that passes through the lattice point. We have 
\[
\inv_\mathrm{initial}(c) = \inv_\mathrm{initial}(c'),\quad \inv_\mathrm{final}(c) = \inv_\mathrm{final}(c') + k-2 a-1,
\]
as the label $k-a$ switches from the last position in $\sigma_{\mathrm{final}}(c')$ to the first one in $\sigma_{\mathrm{final}}(c)$.
\item[D)] Similar to C), but now 
\[
\inv_\mathrm{initial}(c) = \inv_\mathrm{initial}(c') + k-2 a - 1,\quad \inv_\mathrm{final}(c) = \inv_\mathrm{final}(c'),
\]
as the label $k-a$ switches from the last position in $\sigma_{\mathrm{initial}}(c')$ to the first position in $\sigma_{\mathrm{initial}}(c)$.
\item[BE)] There is no difference between $c$ and $c''$. In passing from $c'$ to $c$ the rightmost label disappears from both $\sigma_{\mathrm{initial}}(c')$, $\sigma_{\mathrm{final}}(c')$. The rightmost label in $\sigma_{\mathrm{initial}}(c')$ is $a+1$, while in $\sigma_{\mathrm{final}}(c')$ it is $a+2$, where $a$ is the number of parts of $c$ to the left of the part which passes through the lattice point. Therefore
\[
\inv_\mathrm{initial}(c) = \inv_\mathrm{initial}(c') - k + a,\quad \inv_\mathrm{final}(c) = \inv_\mathrm{final}(c') - k + a+1.
\]
\end{itemize}

Now we see that both sides of the statement satisfy the same recursions and the same initial conditions, so the proof is complete.
\end{proof}

\section{Proof of the compositional $(m,n)$-shuffle conjecture}
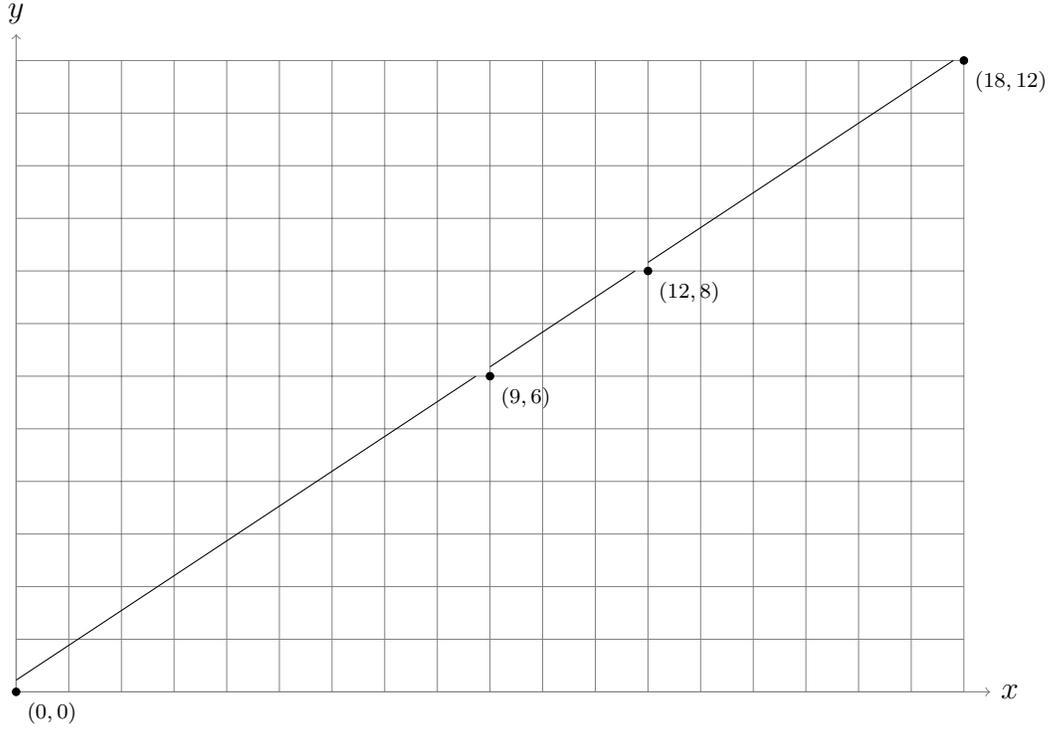
\begin{figure}
\def\eps{0.2}
\def\s{(2/3-0.005)}
\def\t{(1/(1+\s))}
\def\m{18}
\def\n{12}
\def\h{(\n-(\m-\eps)*\s)}
\begin{tikzpicture}[scale=0.7]
\draw[gray, step=1, ultra thin] (0,0) grid (\m,\n);
\draw[gray,thin,->] (0, 0) -- (\m+0.5, 0) node[black,right] {$x$};
\draw[gray,thin,->] (0, 0) -- (0, \n+0.5) node[black,above] {$y$};
\draw (0, {\h}) -- ({(6-\h) / \s},6);
\draw (9, {\h+9*\s}) -- ({(8-\h) / \s},8);
\draw (12, {\h+12*\s}) -- ({(12-\h) / \s},12);
\draw (0,0) node[draw,circle, inner sep=1pt,fill] {} node[below right] {\tiny $(0,0)$};
\draw (9,6) node[draw,circle, inner sep=1pt,fill] {} node[below right] {\tiny $(9,6)$};
\draw (12,8) node[draw,circle, inner sep=1pt,fill] {} node[below right] {\tiny $(12,8)$};
\draw (18,12) node[draw,circle, inner sep=1pt,fill] {} node[below right] {\tiny $(18,12)$};
\end{tikzpicture}
\caption{The coloring corresponding to slope $2/3-\varepsilon$, composition $(3,1,2)$.}
\label{fig:composition323}
\end{figure}
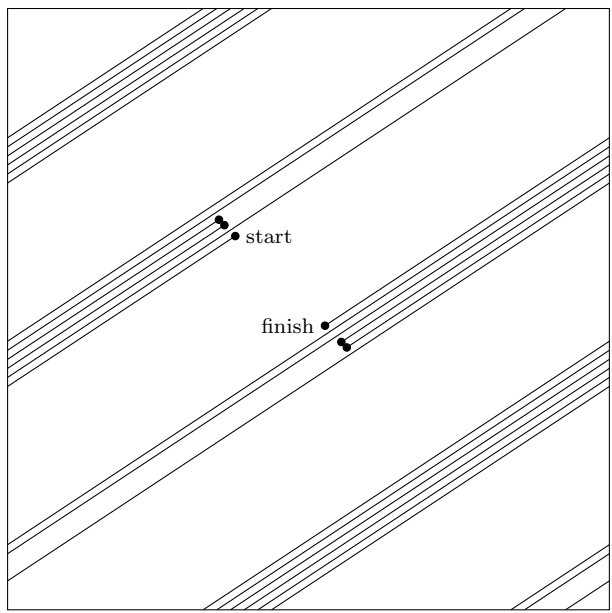
\begin{figure}
\def\eps{0.05}
\def\s{(2/3-0.005)}
\def\t{(1/(1+\s))}
\def\m{18}
\def\n{12}
\def\h{(\n-(\m-\eps)*\s)}
\begin{tikzpicture}[scale=8]
\draw (0,0) rectangle (1,1);

\draw ({(\t)*(1-\h)}, {1-(\t)*(1-\h)}) node[left]{\tiny finish} node[draw,circle, inner sep=1pt,fill] {} -- (1, {\h+\s});
\draw (0, {\h+\s}) -- ({1/\s*(1-\h-\s)}, 1);
\draw ({1/\s*(1-\h)-1)}, 0) -- (1, {\h+2*\s-1});
\draw (0, {\h+2*\s-1}) -- ({1/\s*(2-\h)-2},1);
\draw ({1/\s*(2-\h)-2)},0) -- (1, {\h+3*\s-2});
\def\h{(\n-(\m-\eps)*\s + 3*\s-2)}
\draw (0, {\h}) -- (1, {\h+\s});
\draw (0, {\h+\s}) -- ({1/\s*(1-\h-\s)}, 1);
\draw ({1/\s*(1-\h)-1)}, 0) -- (1, {\h+2*\s-1});
\draw (0, {\h+2*\s-1}) -- ({1/\s*(2-\h)-2},1);
\draw ({1/\s*(2-\h)-2)},0) -- (1, {\h+3*\s-2});
\def\h{(\n-(\m-\eps)*\s + 6*\s-4)}
\draw (0, {\h}) -- (1, {\h+\s});
\draw (0, {\h+\s}) -- ({1/\s*(1-\h-\s)}, 1);
\draw ({1/\s*(1-\h)-1)}, 0) -- (1, {\h+2*\s-1});
\draw (0, {\h+2*\s-1}) -- ({(1/\s*(2-\h)-2)*(1-\t)},{1-(1/\s*(2-\h)-2)*(1-\t)}) node[draw,circle, inner sep=1pt,fill] {};
\def\h{(\n-(\m-\eps)*\s + 9*\s-6)}
\draw ({(\t)*(1-\h)}, {1-(\t)*(1-\h)}) node[draw,circle, inner sep=1pt,fill] {} -- (1, {\h+\s});
\draw (0, {\h+\s}) -- ({1/\s*(1-\h-\s)}, 1);
\draw ({1/\s*(1-\h)-1)}, 0) -- (1, {\h+2*\s-1});
\draw (0, {\h+2*\s-1}) -- ({(1/\s*(2-\h)-2)*(1-\t)},{1-(1/\s*(2-\h)-2)*(1-\t)}) node[draw,circle, inner sep=1pt,fill] {};
\def\h{(\n-(\m-\eps)*\s + 12*\s-8)}
\draw ({(\t)*(1-\h)}, {1-(\t)*(1-\h)}) node[draw,circle, inner sep=1pt,fill] {} -- (1, {\h+\s});
\draw (0, {\h+\s}) -- ({1/\s*(1-\h-\s)}, 1);
\draw ({1/\s*(1-\h)-1)}, 0) -- (1, {\h+2*\s-1});
\draw (0, {\h+2*\s-1}) -- ({1/\s*(2-\h)-2},1);
\draw ({1/\s*(2-\h)-2)},0) -- (1, {\h+3*\s-2});
\def\h{(\n-(\m-\eps)*\s + 15*\s-10)}
\draw (0, {\h}) -- (1, {\h+\s});
\draw (0, {\h+\s}) -- ({1/\s*(1-\h-\s)}, 1);
\draw ({1/\s*(1-\h)-1)}, 0) -- (1, {\h+2*\s-1});
\draw (0, {\h+2*\s-1}) -- ({(1/\s*(2-\h)-2)*(1-\t)},{1-(1/\s*(2-\h)-2)*(1-\t)}) node[draw,circle, inner sep=1pt,fill] {} node[right]{\tiny start};

\end{tikzpicture}
\caption{The braid corresponding to the coloring on Figure \ref{fig:composition323}}
\label{fig:composition323braid}
\end{figure}
Next we specify Theorem \ref{thm:main} to the case of ``colorings with minimal gaps'', as in Remark \ref{rem:compositioncase}. See Figures \ref{fig:composition323} and \ref{fig:composition323braid} for an example, which corresponds to the $(18,12)$-shuffle conjecture with composition $(3,1,2)$.

Let $m,n\in\Z_{>0}$ be relatively prime. Recall the constants $s=\frac{n}{m}-\varepsilon$, $t=\frac{1}{s+1}$. For a composition $\alpha$ we set $B_{n/m}^\alpha:=B_{s,c_\alpha}$,
where $c_\alpha$ was defined in Remark \ref{rem:compositioncase}. In the initial position all points are clustered close to the ``start'' at point $(1-t,t)$, and their labels are in the opposite order. In the final position all points are clustered close to ``finish'' at $(t, 1-t)$, and their labels are again in the opposite order. In particular, the $q$-power factor in Theorem \ref{thm:main} can be ignored. 

Let $b_{m,n}\in \B_1^+(\TT_0)$ denote the braid consisting of a single strand, which starts at $(1-t-\varepsilon, t+\varepsilon)$, travels at slope $s$ downwards until it reaches  a position at $(t-\varepsilon', 1-t+\varepsilon')$ for a small $\varepsilon'>0$ having crossed the horizontal wall $n-1$ times and the vertical wall $m-1$ times. The monoid $\B_1^+(\TT_0)$ is freely generated by $y_1$ and $z_1$, so $b_{m,n}$ is just a sequence of $y_1$'s and $z_1$'s, we write $b_{m,n}(y_1, z_1)$. Denote $b_{m,n}^{(1)}=b_{m,n}$ and more generally
\[
b_{m,n}^{(a)} = (b_{m,n} y_1 z_1)^{a-1} b_{m,n} \quad (a\geq 0).
\]
Now $b_{m,n}^{(a)}$ gives a braid with one strand which begins at $(1-t-\varepsilon, t+\varepsilon)$ and crosses the horizontal resp. vertical walls $an-1$ resp. $am-1$ times, so it is the braid corresponding to the composition with one part $(a)$.

Let $\alpha=(\alpha_1,\alpha_2,\ldots,\alpha_k)$ be a composition and $\alpha'=(\alpha_1,\alpha_2,\ldots,\alpha_k, a)$. Denote the corresponding braids by $B$ and $B'$ respectively. Our aim is to express $B'$ in terms of $B$. First we let the extra new point stay fixed and follow all the steps of $B$. By Proposition \ref{prop:phiplus2} the result is expressed as
\begin{equation}\label{eq:cases1}
\begin{cases}
\varphi_+^*(B)\; T_{1\upto k+1}^* & \text{if $1-t<t$ ( $\mathrm{start}<\mathrm{finish}$),}\\
T_{k+1\downto 1} \;\varphi_+^*(B)\; T_{1\upto k+1}^* & \text{if $t<1-t$ ($\mathrm{finish}<\mathrm{start}$).}
\end{cases}
\end{equation}
Now we move the extra point. Note that each time before we cross a vertical wall the position is $1$, and before we cross a horizontal wall the position is $k+1$. This holds with one exception: each time we are passing next to the finish, the position is $k+1$ instead of $1$. Thus the rules \eqref{eq:braidnext} produce combinations of 
\[
z_1,\; T_{k+1\downto 1},\; \tilde y_{k+1},\; T_{1\upto k+1}^*,
\]
which are arranged in such a way that each sequence of consecutive $\tilde y_{k+1}$ has $T_{1\upto k+1}^*$ on the left and $T_{k+1\downto 1}$ on the right. Note that if the sequence begins with $\tilde y_{k+1}$, then $t<1-t$ so we are in the second case of \eqref{eq:cases1}, and still have an extra $T_{k+1\downto 1}$ on the right. Taking into account the exception, we obtain that when the point is next to the finish for the first time, the braid looks like
\[
T_{k+1\downto 1}\; b_{m,n}(y_1, z_1)\;\varphi_+^*(B)\; T_{1\upto k+1}^*.
\] 
If $a=1$, we are done. If $a>1$, then each extra round begins with
\[
\tilde y_{k+1} z_{k+1} = \tilde y_{k+1}\; T_{k+1\downto 1}\; T_{1\upto k+1}^* \; z_{k+1} 
= \tilde y_{k+1}\; T_{k+1\downto 1}\; z_1\; T_{1\upto k+1},
\]
and continues as before, so the final formula is
\[
B' = \left(T_{k+1\downto 1}\; b_{m,n}(y_1, z_1) y_1 z_1\; T_{1\upto k+1}\right)^{a-1} \;
T_{k+1\downto 1}\; b_{m,n}(y_1, z_1)\;\varphi_+^*(B)\; T_{1\upto k+1}^*.
\]
Applying to $d_+^{k+1}(1)$, in the case $a=1$ we have
\[
B' d_+^{k+1}(1) = T_{k+1\downto 1}\; b_{m,n}(-y_1, (qt)^{-1} z_1)\;\varphi_+^*(B)\; d_+^{k+1}(1)
\]
\[
= T_{k+1\downto 1}\; b_{m,n}(-y_1, (qt)^{-1} z_1)\; (-y_1 d_+^*)\; B d_+^k (1).
\]
It is easy to see that $b_{m,n}$ can be computed recursively using the extended Euclidean algorithm. We have, in particular,
\[
b_{m,n}(-y_1, (qt)^{-1} z_1)\; (-y_1 d_+^*) = (-1)^{m-1} \rho_{m,n}^*(d_+),
\]
\[
b_{m,n}(-y_1, (qt)^{-1} z_1)\; (-y_1 z_1) = (-1)^{m-1} \rho_{m,n}^*(y_1).
\]
Therefore, for any $a>0$
\[
B' d_+^{k+1}(1) = (-1)^{(m-1)a} (qt)^{1-a} \rho_{m,n}^*(y_{k+1}^{a-1} \; T_{k+1\downto 1}^*\; d_+) B d_+^k(1).
\]
By induction, we obtain, for a composition $\alpha=(\alpha_1,\ldots,\alpha_k)$ with sum $N$
\[
D_{n/m, c_\alpha} = B_{n/m}^\alpha d_+^k (1) = (-1)^{(m-1)N} (qt)^{k-N} \rho_{m,n}^*(y_1^{\alpha_1-1} y_2^{\alpha_2-1} \cdots y_k^{\alpha_k-1} d_+^k)(1).
\]
By \eqref{eq:compstocolorings} this gives
\[
D_{n/m, \alpha} = (-1)^{(m-1)N} q^{k-N} \rho_{m,n}^*(y_1^{\alpha_1-1} y_2^{\alpha_2-1} \cdots y_k^{\alpha_k-1} d_+^k)(1),
\]
and matches \eqref{eq:compshufflhs}. Thus the conjecture is proved.

\printbibliography


\end{document}